\documentclass[12pt]{amsart}
\usepackage{amsthm,amsmath,amssymb,amscd,graphics,enumerate}
\input xy
\xyoption{all}

{
   \newtheorem{theorem}[subsubsection]{Theorem}
   \newtheorem{proposition}[subsubsection]{Proposition}     
   \newtheorem{lemma}[subsubsection]{Lemma}
   \newtheorem*{claim}{Claim}
   \newtheorem{corollary}[subsubsection]{Corollary}

   \newtheorem{problem}[subsubsection]{Problem}
}
{\theoremstyle{definition}

   \newtheorem{definition}[subsubsection]{Definition}
   \newtheorem{remark}[subsubsection]{Remark}
}
\newcommand{\RR}{{\mathbb{R}}}

\newcommand{\CC}{{\mathbb{C}}}
\newcommand{\QQ}{{\mathbb{Q}}}

\newcommand{\PP}{{\mathbb{P}}}
\newcommand{\ZZ}{{\mathbb{Z}}}

\newcommand{\bbA}{{\mathbb{A}}}

\newcommand{\bL}{{\mathbf{L}}}

\newcommand{\bR}{{\mathbf{R}}}

\newcommand{\cC}{{\mathcal C}}

\newcommand{\cE}{{\mathcal E}}
\newcommand{\cF}{{\mathcal F}}
\newcommand{\cG}{{\mathcal G}}
\newcommand{\cH}{{\mathcal H}}

\newcommand{\cN}{{\mathcal N}}
\newcommand{\cO}{{\mathcal O}}
\newcommand{\cP}{{\mathcal P}}

\newcommand{\cS}{{\mathcal S}}
\newcommand{\cT}{{\mathcal T}}

\newcommand{\Spec}{\operatorname{Spec}}
\newcommand{\Supp}{\operatorname{Supp}}

\newcommand{\Hom}{{\operatorname{Hom}}}
\newcommand{\cHom}{{{\cH}om}}

\newcommand{\Ker}{{\operatorname{Ker}}}
\newcommand{\Coker}{{\operatorname{Coker}}}
\newcommand{\im}{{\operatorname{im}}}

\newcommand{\Tor}{{\operatorname{Tor}}}

\newcommand{\das}{\dashrightarrow}
\newcommand{\lrar}{\longrightarrow}

\newcommand{\can}{{\operatorname{can}}}

\newcommand{\double}{\genfrac..{0pt}1
{\raise -1pt\hbox{$\scriptstyle\longrightarrow$}}{\raise 3pt\hbox
{$\scriptstyle\longrightarrow$}}} 
\newcommand{\setmin}{\,\protect%
\begin{picture}(8,3.5)\qbezier(1,3.5)(4,2.)(7,.5)\end{picture}\,}

\renewcommand{\setminus}{\setmin}

\def\tototi{\mathbin{\mathop{\otimes}\limits^{\raise-1pt\hbox
{$\scriptscriptstyle {\rm L}$}}}}

\def\indlim{\mathop{\vrule width0pt height7pt depth
4pt\smash{\lim\limits_{\raise 1pt\hbox to 14.5pt
{\rightarrowfill}}}}}
\def\projlim{\mathop{\vrule width0pt height7pt depth
4pt\smash{\lim\limits_{\raise 1pt\hbox to 14.5pt
{\leftarrowfill}}}}}

\newcommand{\ra}{\rightarrow}
\newcommand{\sub}{\subset}
\newcommand{\ov}{\overline}
\newcommand{\id}{\operatorname{id}}
\newcommand{\Om}{\Omega}
\newcommand{\hra}{\hookrightarrow}
\newcommand{\wt}{\widetilde}

\begin{document}

\title[$t$-structures and
  valuative criteria]{Sheaves of $t$-structures and
  valuative criteria for stable complexes}   
\author[Abramovich]{Dan Abramovich}
\thanks{Research of D.A. partially supported by NSF grant DMS-0335501}  
\address{Department of Mathematics, Box 1917, Brown University,
Providence, RI, 02912} 
\email{abrmovic@math.brown.edu}
\author[Polishchuk]{Alexander Polishchuk}
\thanks{Research of A.P. partially supported by NSF grant DMS-0302215}  
\address{Department of Mathematics, University of Oregon, Eugene, OR 97403}  
\email{apolish@uoregon.edu, apolish@math.bu.edu}

\date{\today}
\maketitle
\setcounter{tocdepth}{1}
\tableofcontents

\section{Introduction}

We work with varieties over a field $k$ of characteristic 0, though many arguments
go through without this assumption. The notation $D(X)$ is reserved for the
{\em bounded} derived category of {\em coherent} sheaves on a variety $X$,
namely $D(X):=D^b_c(X)$.

The core of this paper is concerned with the construction of a
``constant'' $t$-structure on the bounded derived category of coherent
sheaves $D(X\times S)$, with $X$ and $S$ smooth varieties, given a
nondegenerate $t$-structure on $D(X)$ with noetherian heart. 

While we believe this construction and the methods involved should be
generally useful in studying derived categories, we present here but one
application: we prove a valuative criterion for separation and
properness for the collection  $\cP(1)$ of stable objects of phase 1 under a 
numerical, locally finite and noetherian
Bridgeland-Douglas stability condition $(Z,\cP)$ on $D(X)$, where
$X$ is a smooth projective variety.  

\subsection{The main results}
Let $X$ be a smooth projective variety and, as always, let $D(X) := D^b_c(X)$
be the 
derived category of bounded complexes of coherent sheaves. Whenever
given a variety $S$ we denote by $p: X\times S \to X$ the projection
to the first factor.

\subsubsection{Results on sheaves of $t$-structures}
 In section
2 we are
given a nondegenerate $t$-structure $(D(X)^{\leq 0}, D(X)^{\geq
0})$ with heart $\cC = D(X)^{\leq 0}\cap D(X)^{\geq
0}$ (see section \ref{Sec:intro-t-structures} for definitions and a brief
introduction). For a smooth projective variety $S$ with ample line bundle $L$
we define  
$$D(X\times S)^{[a,b]}\quad  =\quad  \big\{\ E\, \in D(X\times S)\
\big|\ \bR p_*\, (E \otimes L^n)\ \in\ D(X)^{[a,b]}\quad \forall n\gg
0\  \big\}.$$
Assuming that $\cC$ is noetherian we prove 
\begin{enumerate}
\item\label{It:is-a-t-structure}
 $(D(X\times S)^{\leq 0}, D(X\times S)^{\geq 0})$ is a
nondegenerate $t$-structure on $D(X\times S)$.
\item It extends to a sheaf of $t$-structures over $S$.
\item It is independent of the choice of $L$.
\item The functor $p^*: D(X) \to D(X\times S)$ is $t$-exact.
\item\label{It:closed-exact}
 For every closed immersion of a smooth variety
$i_T:T\hookrightarrow S$ the functor $i_{T\,*}:  D(X\times T)\to
D(X\times S)$ is $t$-exact.
\item\label{It:open-independent}
 If $S'\subset  S$ is open, then the  $t$-structure $(D(X\times
S')^{\leq 0}, D(X\times S')^{\geq 0})$ on  $ D(X\times S')$ is independent of
the 
projective completion $S' \subset S$. 
\item\label{It:heart-noetherian} The heart $\cC_{S'}=  D(X\times
S')^{\leq 0}\cap D(X\times S')^{\geq 0}$ is noetherian. 

And finally
\item\label{It:heart-characterization}
 If $\cC$ is bounded with respect to the standard $t$-structure
on $D(X)$, then we can characterize the $t$-structure by \\
$ E \in D(X\times S)^{\geq 0}\quad \quad \Longleftrightarrow$  
 \begin{center}$ \Hom\big(F, \bR p'_*
(j^*E) \big) = 0 $, whenever $ F \in D(X)^{< 0}$ and $ j: S' \hookrightarrow  S$
 open.
\end{center}
\end{enumerate} 
This is the subject of section \ref{Sec:constant-t}. Parts
(\ref{It:is-a-t-structure})-(\ref{It:closed-exact})  and (\ref{It:heart-noetherian} )  are in  
Theorem \ref{projtstrthm}. Part (\ref{It:heart-characterization}) is
in  Theorem \ref{Th:quasiproj-characterisation}, and
(\ref{It:open-independent}) is in Theorem \ref{Prop:open-independent} (but see
also Theorem \ref{Th:quasiproj-characterisation}).

\subsubsection{Results on flat extensions of objects}

In section \ref{Sec:families-in-heart} we investigate further the 
sheaf of $t$-structures   discussed above. We denote the heart of this
$t$-structure over a smooth quasi-projective scheme $S$ by $\cC_S$.

 Apart from a proof of part
\ref{It:open-independent}  above, we introduce the notion of a {\em family of
objects in $\cC$ parametrized by $S$} (Definition \ref{Def:family-in-C}) and
show a few basic properties:
\begin{enumerate}
\item {\bf The open heart property:} if the restriction of $E\in D(X\times S)$
over a smooth closed subvariety $T\subset S$ is in $\cC_T$ then $E_U\in \cC_U$
for some neighborhood $U$ of $T$ (Proposition \ref{Prop:open-heart}).
\item {\bf The total family is in the heart:} if $E\in D(X\times S)$ is a
family of objects in $\cC$ then $E\in \cC_S$ (Corollary
\ref{Cor:family-is-in-CS}). 
\item {\bf Extension across a divisor:} if $D\subset S$ is a smooth divisor
with complement $U$ and if $E_0\in \cC_U$, then there is an extension $E \in
\cC_S$ of $E_0$ such that $\bL i_D^* E \in \cC_D$ (Proposition
\ref{Prop:extension-across-divisor}).
\end{enumerate}

In section \ref{Sec:valuative-P} we fix a Bridgeland--Douglas stability
condition $(Z,\cP)$ on $D(X)$ (see section \ref{Sec:intro-stability-conditions}
for definitions and a brief introduction).  {\em We assume that its heart is
noetherian}, and  concentrate on a one-parameter 
family $E_0$ of objects in $\cP(1)$ parametrized by an open curve $U\subset S$.
We prove a ``valuative criterion for properness and separation'': 
\begin{enumerate}
\item {\bf valuative properness:} a $t$-flat extension $E$ of $E_0$ over $S$
has its fibers in $\cP(1)$; 
\item {\bf separation up to S-equivalence:} any two extensions have
S-equivalent fibers, and 
\item {\bf polystable replacement:} after a base change there is an extension
with polystable fibers on $S 
\setminus U$. 
\end{enumerate} 
See Theorem \ref{Th:valuative-P1}.

\subsubsection{Examples and questions}
For our results to be of some practical relevance, we show that there is a
large source of examples where they apply. For this purpose we prove in
Proposition 
\ref{Prop:discrete-is-noetherian} that 
\begin{quote} every ``discrete'' stability condition
aligned with $\RR \subset \CC$ has noetherian heart. \end{quote}
In particular, our result
always applies in the case when $X$ is a curve, with a proof very different from the 
classical one. Also, in the case when $X$ is a K3 surface
there exists a connected component in the
space of stability conditions (this component is described in \cite{Bridgeland2}) 
in which
``discrete" stability conditions are dense (see Corollary \ref{rat-cor}).

In section \ref{Sec:questions}, we state a number of natural questions
about moduli spaces, and about possibilities for extending the scope of
validity and reach of our results.

\subsection{$t$-structures}\label{Sec:intro-t-structures}
Let $T$ be a triangulated category.
Recall that a {\it  
$t$-structure} 
$$ (T^{\leq 0}, T^{\geq 0})$$ is determined by a full subcategory
$T^{\leq 0}\subset T$ such that   
\begin{enumerate}
\item $T^{\leq 0}$ is stable under positive (left) shifts,
 i.e. if
  $a\in T^{\leq 0}$ then $a[1]\in T^{\leq 0}$,
and
\item The embedding $T^{\leq 0} \subset T$ has a right adjoint {\em
truncation functor}  
  $\tau_{\leq 0}: T \to T^{\leq 0}$.
\end{enumerate} 
Recall that for a full subcategory $U\subset
T$ the right orthogonal full subcategory is defined by
$$U^\perp \ \ :=\ \  \big\{F \in T\ |\  \Hom(E, F) = 0\quad \forall E
\in U \big\}.$$  Using this 
one defines as usual $T^{\geq 1} = (T^{\leq 0})^\perp$, $T^{\leq -i} =
T^{\leq 0} [i]$ and similarly for $T^{\geq -i}$. 

The heart of the $t$-structure is the full subcategory $ \cC = T^{\leq
  0}\cap T^{\geq 0}.$ 

 The embedding
$T^{\geq 0}\subset T$ has a left adjoint $\tau_{\geq 0}: T \to T^{\geq
0}$, and the composite functor $H^0=\tau_{\geq 0}\tau_{\leq 0}: T \to \cC$ is
the {\em 0-th cohomology functor} corresponding to the
$t$-structure. As usual we denote $H^i(F) = H^0(F[i])$.

The $t$-structure is said to be {\em nondegenerate}\footnote{We use a
  condition stronger than \cite{BBD}. Our definition agrees with the 
  standard conjunction    ``bounded and nondegenerate'', but we use
  ``bounded'' in a slightly   different ``relative'' context below.}
  if  
\begin{enumerate}
\item $\cap_i T^{\leq i} = \cap_i T^{\geq i} = 0$, and
\item $\cup_i T^{\leq i} = \cup_i T^{\geq i} = T$.
\end{enumerate}

There are no definite results about when there is an equivalence $T
\simeq D^b(\cC)$. However, when the $t$-structure is nondegenerate we
have, for $a,b \in \cC$, the expected equalities 
$$\begin{array}{rcll}\Hom^i_{\cC}(a,b) &=&\Hom^i_T(a,b) =  0 &\quad\mbox{ for }i<0,\\
\quad \mbox{and} \quad 
\Hom^i_{\cC}(a,b) &=& \Hom^i_T(a,b)& \quad \mbox{   for }i=0,1. \end{array}$$
Another construction which always works is the following: let 
$F^\bullet$ be a complex of objects in the heart $\cC$ of a $t$-structure on $T$. Then there exists {\em a convolution} $[F^\bullet]\in T$ of $F^\bullet$, unique up to isomorphism (see \cite{Orlov}, Lemma 1.5). It is not difficult to see that the cohomology objects with respect to the $t$-structure satisfy $ H^i([F^\bullet]) = H^i(F^\bullet)$.

We use the standard notation $T^{[a,b]}$ for the full subcategory  $T^{\geq a} \cap T^{\leq b}$. 
Another notion we will use is the following: a class of objects $\cE\subset T$
is said to be {\em bounded with respect to a $t$-structure $(T^{\leq 0},
T^{\geq 0})$ } if there is an integer $N$ such that every object $E \in \cE$
satisfies $$E \in T^{[-N,N]}.$$ We will use this terminology when considering whether the heart $\cC$ of one $t$-structure is bounded with respect to another $t$-structure.

The prototypical example of a $t$-structure is the derived  category $D(Sh(X))$
of sheaves on a 
space, with the usual truncated subcategories 
$\big(\,D^{\leq 0}(Sh(X)),D^{\geq 0}(Sh(X))\,\big)$. This $t$-structure is 
nondegenerate when restricted to the 
bounded derived category $D^b(Sh(X))$. The heart is the abelian category 
of sheaves.

Let $F: T \to S$ be an exact functor of triangulated categories
equipped with  $t$-structures. The functor is said to be
{\em left $t$-exact} (or just {\em left exact}) if it sends $T^{\geq 0}$ to $S^{\geq 0}$, and
analogously for {\em right $t$-exact} and {\em $t$-exact}. 

Given an exact functor of triangulated categories $j^*:D \to D_U$  with
{\em fully faithful} left adjoint $j_!: D_U\to D$ and a right adjoint $j_*: D_U
\to D$, consider the subcategory $$D_F\ \:=\ \ \big\{ a\in D \ |\ j^*(a) =
0\big\}.$$
This notation should remind of the important examples when $D$ is the
category of sheaves over some space, $D_U$ (resp., $D_F$) is the category of sheaves
on an open part (resp., on the complementary closed part).  
A theorem of  Beilinson, Bernstein and Deligne \cite{BBD} says that any
two 
$t$-structures, $(D_F^{\leq 0}, D_F^{\geq 0})$ on $D_F$ and $(D_U^{\leq 0},
D_U^{\geq 0})$ on $D_U$, can be glued together to a $t$-structure on $D$,
such that both the embedding $i_*:D_F\to D$ and the functor $j^*:D\to D_U$ are
$t$-exact. Namely, in this situation the embedding $i_*$ admits left and right adjoint
functors $i^*,i^!:D\to D_F$ and the glued $t$-structure is defined by
$$D^{\leq 0}=\{a\in D\ |\ j^*(a)\in D_U^{\leq 0}, i^*(a)\in D_F^{\leq 0}\},$$ 
$$D^{\geq 0}=\{a\in D\ |\ j^*(a)\in D_U^{\geq 0}, i^!(a)\in D_F^{\geq 0}\}.$$

This has an inductive generalization (again in \cite{BBD}) in the case when a triangulated category $T$ has a semiorthogonal decomposition $$(T_1,T_2,\ldots,T_m).$$
By definition, this means that $T_1,T_2,\ldots,T_m$ are full triangulated subcategories generating
$T$ such that $\Hom(a_i,a_j)=0$ whenever $a_i\in T_i$ and $a_j\in T_j$, where $i>j$.
In addition, we require that for every $n\le m$ 
the embedding into $T$ of the full triangulated subcategory generated by $T_1,T_2,\ldots,T_n$
admits left and right adjoint functors.
In the situation considered above we can view $(D_F,j_!D_U)$ as a semiorthogonal decomposition
of $D$. Applying the gluing construction iteratively 
we can glue $t$-structures $(T_i^{\leq 0},T_i^{\geq 0})$
given on each $T_i$ into a canonical $t$-structure on $T$. 
The glued $t$-structure is characterized by the equality
$$D^{\geq 0}=\{a\in T\ |\ i^!_n(a)\in T_n^{\geq 0}, n=1,\ldots,m\},$$
where $i_n^!:T\to T_n$ is the right adjoint functor to the embedding $T_n\sub T$.

\subsection{Stability conditions}\label{Sec:intro-stability-conditions}
Motivated by work from String Theory 
(see e.g. \cite{Douglas1, Aspinwal-Douglas, Douglas2}), Tom Bridgeland
introduced in \cite{Bridgeland} the a notion of {\em stability condition}
on a triangulated category $T$.

First, a collection of full subcategories $\cP(t), t \in \RR $ is
called a {\em slicing} if 
\begin{enumerate}
\item for all $t\in \RR $ we have $\cP(t+1) = \cP(t)[1]$,
\item if $t_1>t_2$ and $A_i \in \cP(t_i)$, then $\Hom(A_1,A_2)=0$, and
\item every $E\in T$ has a {\em Harder-Narasimhan filtration}, i.e.
there is a finite sequence $t_1> t_2>\cdots>t_n$ and a diagram of distinguished
triangles
$$\xymatrix{0\ar@{=}[r]&E_0 \ar[rr]&& E_1\ar[rr]\ar[dl]&& E_2\ar[r]\ar[dl]&
\cdots\ar[r]& E_{n-1}\ar[rr]&& E_n\ar[dl]\ar@{=}[r]&E\\
&&A_1\ar[lu]^{[1]}&&A_2\ar[lu]^{[1]}&&&&A_n\ar[lu]^{[1]}}$$
with $A_i \in \cP(t_i)$.
\end{enumerate}

The objects of $\cP(t)$ are called {\em semistable of phase $t$}.

Given a slicing $\cP$, one defines, for each interval $I\subset \RR $
(possibly infinite, with either end open or closed), a subcategory
$\cP(I)$ as the minimal extension-closed subcategory containing all
$\cP(t), t\in I$. The category $\cP(\,(0, \infty)\,)\subset T$ determines a
$t$-structure, with  heart $\cP(\,(0,1]\,)$. We call $\cP(\,(0,1]\,)$
the {\em heart} of the slicing $\cP$. 
When $a<b<a+1$ the
category $\cP(\,(a,b)\,)$ is only {\em quasi-abelian} (see
\cite{Bridgeland}, section 4).

A slicing is called {\em locally finite} if there exists  $\eta>0$ such
that for all $t\in \RR $ the quasi-abelian category $\cP(\,(t-\eta, t+\eta)\,)$
is of finite length (with respect to strict short exact sequences, see
\cite{Bridgeland}, Definition 4.1).

A {\em stability condition} $(Z, \cP)$ on $T$ consists of a slicing
$\cP$ and a group homomorphism $Z: K_0(T) \to \CC$ such if $E\in
\cP(\phi)$ we have $$Z(E)\ =\ m(E)\,\cdot\, e^{i\pi\,\phi}$$ with
$m(E) > 0$. The value $m(E)= |Z(E)|$ is called ``the mass of $E$'' in the
literature.  The stability condition is {\em locally finite} if the 
slicing $\cP$ is.

Suppose $T$ is of finite type over a field $k$. Then a stability
condition $(Z,\cP)$ is said to be  {\em numerical} if the function $Z$ factors
through $\cN(T) = K(T)/K(T)^\perp$, where the perpendicular is taken
with respect to the Euler bilinear form $$\chi(E,F)\ \  =\ \  \sum_i\
(-1)^i\,\,\dim_k\,\Hom_T(E,F[i]).$$ 

In \cite{Bridgeland}, Theorem 1.2 and Corollary 1.3 it is shown that,
assuming $\cN(T)$ is finitely generated, 
the collection of 
all locally finite numerical stability conditions forms a complex 
manifold. 

Following Bridgeland's work and inspired by it, the paper \cite{GKR} appeared where more general stability data on triangulated categories are investigated. Our results below apply without change to maximal-slope semistable objects subject to  a noetherian GKR-stability data.


\subsection{Acknowledgments} We are indebted to Tom Bridgeland for
his inspiring work \cite{Bridgeland} and for a number of helpful
conversations. We thank Michael Artin, 
Johan de Jong, Maxim Kontsevich and Amnon Yekutieli for their suggestions
and help.   We finally thank the referee for pointing to a number of unclear points in the submitted version.

\section{Constant families of $t$-structures}\label{Sec:constant-t}

Throughout this section we fix a smooth projective variety
$X$ over $k$ and a nondegenerate $t$-structure $(D^{\le 0},D^{\ge 0})$
on $D(X)$. Let $S$ be a smooth (quasiprojective) variety over $k$.
We always denote by $p:X\times S\ra X$ the natural projection.
Our goal is to construct a $t$-structure
on $D(X\times S)$ 
such that the functor $p^*:D(X)\ra D(X\times S)$ would be $t$-exact, as
well as the restriction  functor $ j^*:D(X\times S)\ra D(X\times U)$ for an
open embedding $j:U \hookrightarrow S$. For a smooth projective variety
$S$ we will construct a canonical such $t$-structure. Given a smooth projective
variety 
$S$ we can thus ``restrict"  such a $t$-structure to $t$-structures
on $D(X\times U)$ for all open subsets $U\subset S$. 

An analogous construction in a  different framework can be found in
\cite{Yekutieli-Zhang}, Theorems 4.11 and 4.14.

\subsection{Sheaves of $t$-structures}\label{sheafsec}

Let $S$ be a smooth variety over $k$. 

\begin{definition}\label{Def:sheaf}
 A {\em sheaf of $t$-structures (on $X$) over $S$} is a collection
of nondegenerate
$t$-structures $$(D(X\times U)^{\le 0},D(X\times U)^{\ge 0})$$
on $D(X\times U)$ for all open subsets $U\subset S$,
such that the restriction functors $D(X\times U)\ra D(X\times U')$
are $t$-exact, where $U'\subset U$.
\end{definition}

\begin{lemma} Assume that we have a sheaf of $t$-structures over $S$
and let $U_i$ be an open covering of $S$. Assume that for
some $F\in D(X\times S)$ we have $F|_{X\times U_i}\in
D(X\times U_i)^{\le 0}$ (resp, $F|_{X\times U_i}\in D(X\times U_i)^{\ge 0}$)
for all $i$. Then $F\in D(X\times S)^{\le 0}$ (resp.,
$F\in D(X\times S)^{\ge 0}$).
\end{lemma}

{\bf Proof.} Note that, given  an object $G\in D(X\times S)$ such that
$G|_{X\times U_i}=0$ for all $i$, we have $G = 0$, since quasi-isomorphisms
are tested on cohomology sheaves with respect to the standard
$t$-structure. 
 By assumption, the restriction functors
commute with the cohomology functors with respect to our $t$-structures, hence
the condition $H^iF=0$ can be checked locally. \qed 

\begin{proposition}\label{tensprop} Assume that we have a sheaf ot
$t$-structures 
over $S$. Then for every vector bundle $V$ on $S$ the
functor $$ \begin{array}{rcl} D(X\times S) & \to & D(X\times S) \\
 F & \mapsto & F \otimes p_S^*V
\end{array}$$
of tensoring with
the pull-back of $V$ to $X\times S$ is $t$-exact.
\end{proposition}

{\bf Proof.} This follows immediately from the previous lemma. \qed

The following theorem implies that
conversely, if the functors of tensoring with pull-backs of line bundles on $S$
are $t$-exact, then the $t$-structure on $D(X\times S)$ extends to a sheaf
of $t$-structures.

\begin{theorem}\label{sheaflinethm}
Let $S$ be a smooth quasiprojective variety, and let
$(D(X\times S)^{\le 0},D(X\times S)^{\ge 0})$ be a nondegenerate
$t$-structure on $D(X\times S)$. Assume that for some ample
line bundle $L$ on $S$ the functor 
\begin{eqnarray*}
D(X\times S) & \longrightarrow & D(X\times S) \\
F & \mapsto &F \otimes p_S^* L
\end{eqnarray*}
 is $t$-exact.
Then there exists a unique extension of the $t$-structure on $D(X\times S)$
to a sheaf of $t$-structures over $S$.
\end{theorem}

We will need several lemmata for the proof.
Let $T\subset S$ be a closed subset.
Denote by $D_T(X\times S)\sub D(X\times S)$ the
full subcategory of objects supported on $X\times T$. This is the same
as the category of objects whose cohomology sheaves are supported in
$X\times T$. 
Note that $D_T(X\times S)$ is a thick  subcategory. The first lemma is
well known:

\begin{lemma}\label{seclem}
Let $f_1,\ldots,f_n$ be sections of some line bundle $M$ on $S$ such that
$T$ is the set of common zeroes of $f_1,\ldots,f_n$.
Then $F\in D_T(X\times S)$
if and only if the morphisms $f_{i_1}\ldots f_{i_d}:F\ra F\otimes M^d$
are zero for all sequences $(i_1,\ldots,i_d)$ of length $d$ for some $d>0$.
\end{lemma}

{\bf Proof.} The ``if" part follows by considering the induced maps
on cohomology sheaves. Let us prove the ``only if" part. 
Using the standard $t$-structure, suppose
$F \in 
D^{[a,b]}$ and apply induction on the length $b-a$, the case $b-a=0$
being immediate. We take $b\leq c<a$ and consider the distinguished
triangle
$$\tau_{\leq c}F \to F \to \tau_{>c} F \to  \tau_{\leq c} F[1],$$ along with 
the associated long exact sequence of $\Hom(F, \bullet)$. By induction
there is an integer $d$ which works for the two truncations. We claim
that $2d$ will work for $F$. Indeed, the arrow 
$$\tau_{>c}\circ (f_{i_1}\ldots f_{i_d}): F \to  \tau_{>c} F\otimes M^d $$
vanishes since it factors through $f_{i_1}\ldots f_{i_d}:\tau_{>c} F
\to  \tau_{>c} F\otimes M^d$. Therefore $f_{i_1}\ldots f_{i_d}: F \to
F\otimes M^d $ 
factors through an arrow $g:F \to \tau_{\leq c}F\otimes M^d$. It follows that
$$(f_{i_{d+1}}\ldots f_{i_{2d}})\circ g:  F \to \tau_{\leq c}F\otimes M^{2d}$$
vanishes, and therefore also its image $f_{i_{1}}\ldots f_{i_{2d}}:
F \to F\otimes M^{2d}$ vanishes as well.\qed

\begin{lemma}\label{suplem}
Under the assumptions of Theorem \ref{sheaflinethm},
let $H^i$, $i\in\ZZ$ be the cohomology functors with respect
to the $t$-structure on $D(X\times S)$.
Then for every closed subset $T\sub S$ one has
$F\in D_T(X\times S)$ if and only if
$H^i F\in D_T(X\times S)$ for all $i$.
\end{lemma}

{\bf Proof.} 
The ``if" part follows as in  Lemma \ref{seclem}  since our
$t$-structure is nondegenerate. 
The ``only if" part follows from Lemma \ref{seclem} and from the
fact that the twisting functor $F\mapsto F\otimes p_S^*L$ is $t$-exact.
Indeed, replacing $L$ by a suitable power we can assume that
$T$ is the set of common zeros of sections $f_1,\ldots,f_n$ of $L$.
Now if the morphism $f=f_{i_1}\ldots f_{i_d}:F\ra F\otimes L^d$ is zero,
then using the defining property of the truncating functor
$\tau_{\le 0}$ (resp., $\tau_{\ge 0}$) we deduce that
$f:\tau_{\le 0}F\ra \tau_{\le 0}F\otimes L^d = \tau_{\le 0}(F\otimes
L^d)$ 
(respecively, 
$f:\tau_{\ge 0}F\ra \tau_{\ge 0}F\otimes L^d = \tau_{\ge 0}(F\otimes
L^d) $) is zero. \qed

\begin{lemma}\label{Serrelem}
Let $\cC\sub D(X\times S)$ be the heart of the $t$-structure on
$D(X\times S)$ considered in Theorem \ref{sheaflinethm}. 
Then for every closed subset $T\sub S$ the
subcategory $D_T(X\times S)\cap\cC\sub\cC$
is closed under subobjects, quotients and extensions.
\end{lemma}

{\bf Proof.} This follows immediately from Lemma \ref{seclem}. \qed

{\bf Proof of Theorem \ref{sheaflinethm}.}
Let $U\sub S$ be an open set with the complement $T$. Proposition 5.5.4
of \cite{TT} implies that the restriction functor 
$D(X\times S)\ra D(X\times U)$ is essentially surjective. Hence, 
the category $D(X\times U)$ is the quotient of
$D(X\times S)$ by the thick subcategory $D_T(X\times S)$.
By the definition, this means that $D(X\times U)$ is the localization
of $D(X\times S)$ with respect to the localizing class consisting
of morphisms $f:F'\ra F$ such that the cone of $f$ belongs to
$D_T(X\times S)$.
To show that the $t$-structure on $D(X\times S)$
induces a $t$-structure on the localized category, it suffices
to show that every diagram
$$F\stackrel{f}{\leftarrow} F'\ra G$$
such that the cone of $f$ is in $D_T(X\times S)$,
$F\in D(X\times S)^{\le 0}$ and $G\in D(X\times S)^{\ge 1}$,
the induced morphism $F\ra G$ in the localized category vanishes.
Considering the long exact sequence of cohomology associated
with the exact triangle $F'\stackrel{f}{\ra} F\ra C\ra F'[1]$ we have that $H^0C \to H^1F'$ is surjective and $H^iC \to H^{i+1}F'$ an isomorphism for $i\geq 1$. 
Using Lemmata \ref{Serrelem} and \ref{suplem} we derive that
$\tau_{\ge 1}(F')\in D_T(X\times S)$. Therefore,
the arrow $\tau_{\le 0}(F')\ra F'$ belongs to our localizing class
and the above diagram is equivalent to
$$F\leftarrow\tau_{\le 0}(F')\ra G.$$
But $\Hom(\tau_{\le 0}(F'),G)=0$ which proves our claim. \qed

For later use we record the following related  (and probably well
known) lemma: 
\begin{lemma}\label{Lem:extend-morphism} Let $U \subset V$ be an open
subset whose complement is a divisor $T$ with defining section $f \in
H^0(V, \cO_V(T))$. Let $F_1$ and $F_2$ be objects in $D(V)$, and let
$\phi_U: (F_1)_U \to ( F_2)_U$ be a morphism. Then for some $k$,
the morphism $f^k \phi_U$ extends to $\phi': F_1 \to F_2$.
\end{lemma}

{\bf Proof.} The morphism $\phi_U$ corresponds to a diagram 
$$ F_1 \leftarrow H \rightarrow F_2$$
where $$C = Cone( H \to  F_1 ) \in D_T(V).$$ In particular there is an
integer $k$ such that the morphism $f^k : C \to C(kT)$ is zero.

In the long exact sequence 
$$\cdots \to \Hom( F_1 , H(kT)) \to \Hom( F_1 , F_1(kT)) \to \Hom( F_1
, C(kT)) \to \cdots$$
the image of the element $f^k \in \Hom( F_1 , F_1(kT))$ in $\Hom( F_1
, C(kT)) $ factors through the zero morphism  $f^k : C \to C(kT)$, therefore 
it
vanishes. Thus there is $\psi\in   \Hom( F_1 , H(kT))$ mapping to
$f^k \in \Hom( F_1 , F_1(kT))$, representing the inverse of $H_U \to
(F_1)_U$, as required.\qed

The next theorem follows from the
proof of Theorem 3.2.4 of \cite{BBD} (we only need
the case of Zariski sheaves of $\cO$-modules,
so we leave to the reader to formulate
a more general version on arbitrary site similar to the
formulation in \cite{BBD}).

\begin{theorem}\label{BBDthm} Let $U_i$ be a finite open covering of $S$,
and let 
$$\big(\,K_i\in D(X\times U_i)\, ,\ \
\alpha_{ij}:K_i|_{X\times U_{ij}}\ \widetilde{\rightarrow}\
K_j|_{X\times U_{ij}}\,\big)$$
 be a gluing datum, where
$U_{i_1\ldots i_n}:=U_{i_1}\cap\ldots\cap U_{i_n}$ (so
$\alpha_{jk}\alpha_{ij}=\alpha_{ik}$ over $X\times U_{ijk}$).
Assume that $\Hom^j(K_i|_U,K_i|_U)=0$ for all $j<0$ and
all open sets of the form $U=U_{ii_1\ldots i_n}$. Then
there exists $K\in D(X\times S)$ 
equipped with isomorphisms $K|_{X\times U_i}\simeq K_i$ compatible
with $(\alpha_{ij})$ on double intersections.
\end{theorem}

For every $F,G\in D(X\times S)$ let us define
an object $\bR \cHom_S(F,G)\in D(S)$ by setting
$$\bR \cHom_S(F,G)\quad =\quad\bR p_{S*}\ \bR \cHom(F,G).$$
Note that for every open subset $U\subset S$ we have
$$\bR \Hom(F_U,G_U)\quad\simeq\quad \bR \Gamma(U\,,\ \bR \cHom_S(F,G)|_U\,),$$
where for $F\in D(X\times S)$ we set $F_U:=F|_{X\times U}$.

\begin{lemma}\label{sheaflem} 
Let $F$ and $G$ be a pair of objects in $D(X\times S)$ 
such that $$\bR \cHom_S(F,G)\ \in\ D(S)^{\ge 0}$$ (with respect to the  
standard $t$-structure). Then $U\mapsto \Hom(F_U,G_U)$ is a sheaf 
on $S$.
\end{lemma}

{\bf Proof.} Indeed, for an open set $U\subset S$ we have
$$\begin{array}{rcl}\Hom(F_U,G_U)\ =\ H^0\,\bR \Hom(F_U,G_U)&=&H^0\,\bR
\Gamma\big(U\,,\,\bR \cHom_S(F,G)|_U\,\big)\\
&=&  
\Gamma\big(U\,,\,H^0\,\bR \cHom_S(F,G)\,\big),\end{array}$$ 
since the functor $\Gamma$ is left exact. \qed

\begin{corollary}\label{sheafcor}  
Let $U_i$ be a finite open covering of $S$, and let 
$$\big (\,K_i\in D(X\times U_i)\,, \ \  
\alpha_{ij}:K_i|_{X\times U_{ij}}\ \widetilde{\rightarrow}\ 
K_j|_{X\times U_{ij}})$$ be a gluing datum.  
Assume that for every $i$ one has
$\bR \cHom_{U_i}(K_i,K_i)\in D(U_i)^{\ge 0}$ (with respect
to the standard $t$-structure). Then
there exists a unique $K\in D(X\times S)$ (up to unique isomorphism)
equipped with isomorphisms $K|_{X\times U_i}\simeq K_i$ compatible with
$(\alpha_{ij})$ on double intersections.
\end{corollary}

{\bf Proof.} Uniqueness follows immediately from Lemma \ref{sheaflem}
since our assumptions imply that $\bR \cHom_S(K,K)\in D(S)^{\ge 0}$.
The existence follows from Theorem \ref{BBDthm}. Indeed, if
$U=U_{ii_1\ldots i_k}$ then the condition
$\bR \cHom_U(K_{i U}, K_{i U})\in D(U)^{\ge 0}$ implies that
$\Hom^j(K_{i U}, K_{i U})=0$ for $j<0$. \qed

\begin{corollary}\label{stackcor}
Let $U\mapsto (D(X\times U)^{\le 0},D(X\times U)^{\ge 0})$
be a sheaf of $t$-structures over $S$. Then the hearts
$\cC_U=(D(X\times U)^{\le 0}\cap D(X\times U)^{\ge 0})$
form a stack of abelian categories on $S$.
\end{corollary}

{\bf Proof.}
Assume that $F,G\in\cC_S$. Then
for every open affine subset $U\subset S$ the vanishing of
$H^i \bR \Hom(F_U,G_U)$ for $i<0$ implies that
$\bR \cHom_S(F,G)|_U$ belongs to $D(U)^{\ge 0}$ (with respect
to the standard $t$-structure on $D(U)$). Hence,
$\bR \cHom_S(F,G)\in D(S)^{\ge 0}$, so by Lemma \ref{sheaflem}
the groups $U\mapsto \Hom(F_U,G_U)$ form a sheaf on $S$.
To glue objects $F_i\in\cC_{U_i}$
given on some finite open covering $(U_i)$ of $S$ and equipped
with the gluing data we apply Corollary \ref{sheafcor}. \qed

\subsection{Categorification of Hilbert's basis theorem}

In this subsection we will prove a technical result that
can be considered as a categorical version of Hilbert's basis
theorem (in its graded version). An analogous result can be found in
\cite{Artin-Zhang} in their context (it may be that their result implies ours).

Let $\cC$ be an abelian $k$-linear category.

\begin{definition} Let $A=\oplus_n A_n$ be a $\ZZ$-graded associative algebra
with unit over $k$ such that $\dim_k A_n<\infty$ for every $n\in\ZZ$.

\begin{enumerate}
\item[(a)] A {\it graded $A$-module in} $\cC$ is a collection
$M=(M_n |\ n\in\ZZ)$ of objects in $\cC$ and a collection of morphisms
$A_m\otimes M_n\ra M_{m+n}$  for
$m,n\in\ZZ$ satisfying the 
natural associativity condition and such that the composition
$k\otimes M_n\ra A_0\otimes M_n\ra M_n$ is the identity morphisms
for every $n\in\ZZ$. Graded $A$-modules in $\cC$ form an abelian category
in a natural way.

\item[(b)] A {\it free} graded $A$-module {\em of finite type} in $\cC$ is
a finite direct sum of graded $A$-modules in $\cC$ of the form
$A\otimes F(i)$, where 
$$(A\otimes F(i))_n=A_{n+i}\otimes F,$$ 
$F$ is an object of $\cC$, $i$ is a fixed integer; the morphisms
$$A_m\otimes (A\otimes F(i))_n\ra (A\otimes F)_{m+n}$$
are induced by the multiplication in $A$.

\item[(c)] A graded $A$-module $M$ in $\cC$ is called of {\it finite type}
if there exists a surjection $P\ra M$ where $P$ is a free graded $A$-module
of finite type in $\cC$. Note that if $A$ is generated by $A_1$ over
$A_0$ and $M_n=0$ for $n<0$ then $M$ is of finite type 
if and only if
the maps $A_1\otimes M_n\ra M_{n+1}$ are surjective for $n\gg 0$.  
\end{enumerate}
\end{definition}

For example, if $\cC$ is the category of finitely generated modules over
a noetherian ring $R$ then a graded $A$-module in $\cC$ is the same as a
graded $R\otimes_k A$-module with graded components finitely generated
over $R$. A graded $A$-module in $\cC$ is of finite type if and only if
the corresponding $R\otimes_k A$-module is finitely generated in the
usual sense.

In the case $A=k[x_1,\ldots,x_d]$, we have the following analogue
of Hilbert's basis theorem (the proof is also completely analogous).

\begin{theorem}\label{catHilbthm} 
Assume that $\cC$ is noetherian and let $A=k[x_1,\ldots,x_d]$
be the algebra of polynomials in $d$ variables with the grading 
$\deg(x_i)=1$. 

\noindent
(i) Let $F$ be an object of $\cC$, $A\otimes F$ be the corresponding
free graded $A$-module in $\cC$. Then every graded submodule in $A\otimes F$
is of finite type.

\noindent
(ii) The category of graded $A$-modules of finite type
in $\cC$ is abelian and noetherian.
\end{theorem}

{\bf Proof.} 
We  use induction in $d$. When $d=0$ both (i) and (ii) are clearly
true. Let $d>0$ and assume that the assertion is true for $d-1$.
Let $B=k[x_2,\ldots,x_d]$, so that $A=B[x_1]$. 

(i) Let $M\sub A\otimes F$ be a graded submodule. This means
that for every $n\ge 0$ we have a subobject $M_n\sub A_n\otimes F$
such that for every $i=1,\ldots,d$ the image of $M_n$ under the
map $A_n\otimes F\stackrel{x_i}{\ra} A_{n+1}\otimes F$ induced by
the multiplication by $x_i$, is contained in $M_{n+1}$. 

First, for every $i$
we are going to define a graded $B$-submodule $L^i\sub B\otimes F$
of ``leading terms'' of $M$ with $x_1^i$. Consider the filtration
$F_0A\sub F_1A\sub\ldots$ of $A$ induced by the degree in $x_1$, so
that $F_iA$ are polynomials of degree $\le i$ in $x_1$. Let
$F_0A_n\sub F_1A_n\sub\ldots$ be the induced filtration on $A_n$.
Note that for every $i$ and $n$
we have a natural isomorphism $F_iA_n/F_{i-1}A_{n-1}\simeq B_{n-i}$
(taking the coefficient with $x_1^i$). Let
$p_{n,i}:F_iA_n\otimes F\ra B_{n-i}\otimes F$ be the induced projection.
Let us set 
$$F_iM_n:=M_n\cap (F_iA_n\otimes F)\sub M_n,$$
$$L^i_n:=p_{n+i,i}(F_iM_{n+i})\sub B_n\otimes F.$$
Since $M$ is closed under multiplication by $B$, the collection
$(L^i_n |\ n\ge 0)$ is a graded $B$-submodule $L^i$ in $B\otimes F$.
Furthermore, since 
$M$ is closed under multiplication by $x_1$, we obtain the inclusion
 $L^i_n\sub L^{i+1}_n\sub B_n\otimes F$ for every $n$.
Thus, we have an increasing chain of graded $B$-submodules
$L^0\sub L^1\sub\ldots$ in $B\otimes F$.
Since the category of graded $B$-modules of finite type
in $\cC$ is noetherian, this chain necessarily stabilizes.
Hence, there exists $N_1>0$ such that $L^i=L^{i+1}$ for $i\ge N_1$.
By induction assumption the $B$-modules $L^0,\ldots,L^{N_1}$ 
are of finite type. Hence, we can find $N_2>0$ such that the morphisms
$B_1\otimes L^i_n\ra L^i_{n+1}$ are surjective for $n\ge N_2$ and 
$i\le N_1$.

We claim that in this case the moprhisms
$\mu_n:A_1\otimes M_n\ra M_{n+1}$ are surjective for $n\ge N=N_1+N_2$.
It suffices to prove that the induced maps
$$\mu_n^i:\ 
\big(\,(x_1)\,\otimes\  F_iM_n/F_{i-1}M_n\,\big) 
\  \oplus  \  \big(\,B_1\,\otimes\
F_{i+1}M_n/F_iM_n\,\big) \quad \lrar  \quad
F_{i+1}M_{n+1}/F_iM_{n+1}$$
are surjective for $0\le i\le n$, where $(x_1)\sub A_1$ is the line spanned
by $x_1$.
Note that since the kernel of $p_{n,i}$ is $F_{i-1}A_n\otimes F$,
we have a natural isomorphism
$$F_iM_n/F_{i-1}M_n\simeq L^i_{n-i}.$$
Moreover, the map $\mu_n^i$ can be identified with the natural map
$$L^i_{n-i} \ \oplus\  \big(\,  B_1\otimes L^{i+1}_{n-i-1}\,\big) \quad 
\lrar \quad L^{i+1}_{n-i}.$$
If $i\ge N_1$ then  $L^i=L^{i+1}$, hence $\mu_n^i$ is surjective.
If $i<N_1$ then $n-i>N_2$ and the result follows from surjectivity
of the map $B_1\otimes L^{i+1}_{n-i-1}\ \ra\  L^{i+1}_{n-i}$.

\noindent
(ii) First, let us check that for every object $F\in\cC$ the corresponding
free graded $A$-module $A\otimes F$ is a noetherian object.
Let $M^1\sub M^2\sub\ldots $ be an increasing chain of graded $A$-submodules
in $A\otimes F$. For every $n\in\ZZ$ the corresponding chain
$M^1_n\sub M^2_n\sub\ldots $ of subobjects in $A_n\otimes F$ stabilizes
since $\cC$ is noetherian. Hence, we can set $M_n:=\cup_i M^i_n$.
Clearly, $M=(M_n)$ is a graded $A$-submodule of $A\otimes F$. According
to (i) it is of finite type. This easily implies that $M=M^i$ for some $i$.
It follows that every free graded $A$-module of finite type is 
noetherian. Therefore, every graded $A$-module of finite type is
noetherian. Conversely, a noetherian graded $A$-module is of finite type
(otherwise, one would get an infinite strictly increasing chain of submodules).
\qed

\subsection{Constant $t$-structures for projective spaces}
We return to the situation where we are given a nondegenerate $t$-structure $(D(X)^{\leq 0}, D(X)^{\geq 0})$ on the bounded derived category $D(X)$ of coherent sheaves on $X$. In this section we construct, under the assumption that the heart is noetherian,  a sheaf of $t$-structures on $X$ over $\PP^r$, as the limit of a sequence of $t$-structures glued by means of natural semiorthogonal decompositions.

\paragraph{\bf A.  Glued $t$-structures}

\begin{proposition} The following is a nondegenerate $t$-structure on $D(X\times\PP^r)$:
$$D^{[a,b]}\ \ =\ \ \big\{\, F\ \big|\ 
\bR p_*(F)\in D^{[a,b]},\ \bR p_*(F(1))\in D^{[a,b]},\
\ldots\,,\ \bR p_*(F(r))\in D^{[a,b]}\ \big\}. $$
\end{proposition}

{\bf Proof.} 
We claim that this $t$-structure is obtained
by gluing from the standard $t$-structures from the semi\-ortho\-gonal
decomposition of $D(X\times\PP^r)$ into subcategories
$$\big (\,p^*D(X)(-r)\,,\ \ldots\,,\ p^*D(X)(-1)\,,\
p^*D(X)\,\big).$$
Indeed, let us consider the full embeddings
$$i_n:D(X)\to D(X\times\PP^r):F\mapsto p^*F(-n),$$
$n=0,\ldots,r$ and let $i_n^!$ be the right 
adjoint functor to $i_n$, so that
$$i_n^!(F)=\bR p_*(F(n)).$$
For every $F$ let us define a sequence
of objects $F_0=F,F_1,\ldots,F_r$ from the exact triangles
\begin{equation}\label{Fn-tr}
i_ni_n^!F_n\to F_n\to F_{n+1}\to i_ni_n^!F_n[1]\ldots,
\end{equation}
where $n=0,\ldots,r-1$. Note that $F_n$ belongs to the triangulated subcategory 
$D[-r,-n]\sub D(X\times\PP^r)$
generated by $(p^*D(X)(-r),\ldots, p^*D(X)(-n))$. Moreover $F_{n+1}=s_n^*F_n$,
where $s_n^*:D[-r,-n]\to D[-r,-n-1]$ is the left adjoint functor to the embedding
$s_n:D[-r,-n-1]\to D[-r,-n]$.  From this one can 
see that the glued $t$-structure is given by the following formulae:
$$D^{\ge 0}=\{ F\ |\ i_0^!F\in D^{\ge 0}, i_1^!F\in D^{\ge 0},\ldots, i_r^!F\in D^{\ge 0}\},$$
$$D^{\le 0}=\{ F\ |\ i_0^!F_0\in D^{\le 0}, i_1^!F_1\in D^{\le 0},\ldots, i_r^!F_r\in D^{\le 0}\}.$$
We claim that in fact
\begin{equation}\label{glued-tstr}
D^{\le 0}=\{ F\ |\ i_0^!F\in D^{\le 0}, i_1^!F\in D^{\le 0},\ldots, i_r^!F\in D^{\le 0}\}.
\end{equation}
Indeed, applying functors $i_m^!$ (for $m\ge n+1$) to the triangle \eqref{Fn-tr} we get the exact triangles
$$H^0(\PP^r,\cO_{\PP^r}(m-n))\otimes i_n^!F_n\to i_m^!F_n\to i_m^!F_{n+1}\to 
H^0(\PP^r,\cO_{\PP^r}(m-n))\otimes i_n^!F_n[1]\ldots,$$
from which \eqref{glued-tstr} can be easily derived.
\qed    

For every $n\in\ZZ$ set

\begin{equation}\label{ntstreq}
{D^{[a,b]}(X\times\PP^r)_n   }\quad  = 
\quad \left\{\ F \ \left|\begin{array}{rcc} \bR p_*(F(n))&\in& D^{[a,b]},\\  \bR p_*(F(n+1)&\in& D^{[a,b]},\\
&\vdots&\\ \bR p_*(F(n+r))&\in& D^{[a,b]} \end{array} \right.\right\}. 
\end{equation}
Note that we have 
$F\in D^{[a,b]}(X\times\PP^r)_n$ if and only if $F(n)\in D^{[a,b]}(X\times\PP^r)_0$.
Hence, by the previous proposition 
$(D^{\le 0}(X\times\PP^r)_n,D^{\ge 0}(X\times\PP^r)_n)$ is a $t$-structure
on $D(X\times\PP^r)$. 
We will denote by $\tau_{\ge i}^n$, $\tau_{\le i}^n$ the truncation
functors with respect to the $n$-th $t$-structure. 
It is clear that for every $n\ge 0$ the functor
$p^*:D(X)\ra D(X\times\PP^r)$ is $t$-exact with respect to the $n$-th
$t$-structure.

\begin{lemma}\label{inclem}
For $m<n$ one has the following inclusions
$$D^{\le 0}(X\times\PP^r)_m\sub D^{\le 0}(X\times\PP^r)_n,$$
$$D^{\ge 0}(X\times\PP^r)_n\sub D^{\ge 0}(X\times\PP^r)_m\sub 
D^{\ge -r}(X\times\PP^r)_n.$$
\end{lemma}

{\bf Proof.} 
We claim that for $m\le n$ there exists an exact sequence of the form 
\begin{equation}\label{prexactseq}
0\ra V_0\otimes\cO_{\PP^r}(m)\ra V_1\otimes\cO_{\PP^r}(m+1)\ra\ldots\ra
V_r\otimes\cO_{\PP^r}(m+r)\ra \cO_{\PP^r}(n+r)\ra 0,
\end{equation}
where $V_r$ are finite-dimensional $k$-vector spaces. 
Indeed, 
consider the polynomial algebra $S=\oplus_{n\ge 0}H^0(\PP^r,\cO_{\PP^r}(n))$.
Then the standard Koszul complex shows that this algebra is {\it Koszul},
i.e., $\Tor^S_n(k,k)$ is concentrated in internal degree $n$ for all $n$ (see \cite{P}). For an $S$-module $M$ graded in degrees $\geq 0$ it follows that the graded module  $\Tor^S_n(k,M)$ lies in internal degrees $\geq n$; such
 a module $M$ is called {\it Koszul} if $\Tor^S_n(k,M)$ is
concentrated precisely in internal degree $n$  (see \cite{BGS}). Equivalently, $M$ should have a linear free resolution.
Let us denote by $M(i)$ the module obtained from $M$
by the degree shift: $M(i)_n=M_{i+n}$. Set also $M_{\ge j}:=\oplus_{n\ge j} M_n$.
It is easy to see that for every Koszul $S$-module $M$ the module $M(1)_{\ge 0}$ is
also Koszul. Indeed, the long exact sequence of $\Tor$-spaces
associated with the exact sequence 
$$0\ra M(1)_{\ge 0}\ra M\ra M_0\otimes k\ra 0$$ gives that   $\Tor^S_n(k,M)$ lies in degrees $n-1,n$, and since $M(1)_{\ge 0}$ is graded in degrees $\geq 0$, the component of degree $n-1$ in $\Tor^S_n(k,M)$ vanishes.
Iterating this operation we obtain that 
for a Koszul $S$-module $M$ all the modules $M(i)_{\ge 0}$ for $i\ge 0$ are still Koszul.
Applying this to $M=S$ we derive that for every $i\ge 0$ the module $S(i)_{\ge 0}$
admits a linear free resolution. Since $S(i)_{\ge 0}$ has projective dimension $\le r$ such a resolution
will have the following form
$$0\ra V_0\otimes S(-r)\ra V_1\otimes S(-r+1)\ra\ldots\ra V_{r-1}\otimes S(-1)\ra
V_r\otimes S\ra S(i)_{\ge 0}\ra 0.$$
Using the localization functor from finitely generated graded $S$-modules to coherent
sheaves on $\PP^r$ and making an appropriate twist we get the required exact sequence.

Let us first consider the sequence \eqref{prexactseq} for $n=m+1$. In this case 
it coincides with the twist of the Koszul resolution associated with the basis of global
sections of $\cO_{\PP^r}(1)$, so that
$$V_i\simeq \bigwedge^{r+1-i} H^0(\PP^r,\cO_{\PP^r}(1)).$$ 
In particular, we have $V_0\simeq k$. Therefore, we can view our exact sequence as a resolution
for $\cO_{\PP^r}(m)$ in terms of 
$\cO_{\PP^r}(m+1),\ldots,\cO_{\PP^r}(m+r+1)$ positioned in degrees $[0,\ldots,r]$.
Thus, for $F\in D^{\ge 0}(X\times\PP^r)_{m+1}$ we conclude that
$\bR p_*(F(m))\in D^{\ge 0}$, hence $D^{\ge 0}(X\times\PP^r)_{m+1}\sub
D^{\ge 0}(X\times\PP^r)_m$. Iterating we get one of the required inclusions
$$D^{\ge 0}(X\times\PP^r)_n\sub D^{\ge 0}(X\times\PP^r)_m$$
for $m<n$. On the other hand, for arbitrary $m<n$
we can view the exact sequence \eqref{prexactseq} as a resolution for $\cO_{\PP^r}(n+r)$
in terms of $\cO_{\PP^r}(m),\ldots,\cO_{\PP^r}(m+r)$ positioned in
degrees $[-r,\ldots,0]$.
Therefore, for $F\in D^{\ge 0}(X\times\PP^r)_m$ we obtain
$\bR p_*(F(n+r))\in D^{\ge -r}$ for every $n>m$.
This implies that 
$$F\in D^{\ge -r}(X\times\PP^r)_{n+r}\sub D^{\ge -r}(X\times\PP^r)_n.$$
The two other inclusions follow by passing to left orthogonals. 
\qed

\paragraph{\bf B.  Stabilization of cohomology of the glued $t$-structures}
  
Assume that we have two $t$-structures, $(D_1^{\le 0}, D_1^{\ge 0})$ and
$(D_2^{\le 0}, D_2^{\ge 0})$, on a triangulated category $D$, such
that $D_1^{\le 0}\sub D_2^{\le 0}$ (and hence
$D_2^{\ge 0}\sub D_1^{\ge 0}$). Let us denote by
$\tau^1_{\le n}$, $\tau^1_{\ge n}$ (respectively, $\tau^2_{\le n}$,
$\tau^2_{\ge n}$) 
the truncation functors with respect to the first (respectively,
second) 
$t$-structure. Then for every object $F\in D$ 
the canonical morphism $\tau^1_{\le 0}(F)\ra F$ factors uniquely through
$\tau^2_{\le 0}(F)$, so we get a morphism
of functors $\tau^1_{\le 0}\ra \tau^2_{\le 0}$. Similarly,
the canonical morphism $F\ra\tau^2_{\ge 0}(F)$ factors through
$\tau^1_{\ge 0}(F)$, so we get a morphism of functors
$\tau^1_{\ge 0}\ra \tau^2_{\ge 0}$. Let $H^i_1$ (resp. $H^i_2$)
be the $i$-the cohomology functor with respect to the first
(resp. second) $t$-structure. Then we get a morphism of functors
$$H^i_1(F)=\tau^1_{\ge 0}\tau^1_{\le 0}(F[i])\ra
\tau^2_{\ge 0}\tau^2_{\le 0}(F[i])=H^i_2(F).$$

\begin{proposition}\label{stabprop} Assume that the heart $D(X)^{\leq 0}\cap
  D(X)^{\geq 0}$ is 
  noetherian. 

Let $H^i_n$, $i\in\ZZ$ be the cohomology functors
associated with the $n$-th $t$-structure defined by (\ref{ntstreq}).
For every $F\in D(X\times\PP^r)$ there exists an integer $N$ such that
for $n>N$ the morphisms $H^i_n(F)\ra H^i_{n+1}(F)$ are isomorphisms
for all $i$.
\end{proposition}

We need some preparation before giving a proof. Let $\cC$ denote the heart of
the original $t$-structure on $X$. Let us also consider the truncated symmetric algebra
$$A_V:=\oplus_0^r S^i V,$$
where $V=H^0(\PP^r,\cO_{\PP^r}(1))$. 
We denote by $\cC_V$ the full subcategory
in the category of graded representations of the algebra $A_V$ in 
$\cC$ formed by representations $M=(M_n)$ with $M_n=0$ for $n\not\in[0,r]$. 
We claim that exact functor 
$$\Phi:D^{[0]}(X\times\PP^r)_0\to\cC_V:F\mapsto \oplus_{i=0}^r \bR p_*(F(i)).$$ 
induces an equivalence of the heart of the $0$-th $t$-structure $D^{[0]}(X\times\PP^r)_0$ 
with $\cC_V$. In fact, we can also construct a functor in the opposite direction 
$\Psi:\cC_V\to D^{[0]}(X\times\PP^r)_0$
as follows:
$$\Psi(M)=[M_0\boxtimes \Om^r_{\PP^r}\to M_1\boxtimes \Om^{r-1}_{\PP^r}(-1)\to
\ldots\to M_{r-1}\boxtimes \Om^1_{\PP^r}(-r+1)\to M_r\boxtimes\cO_{\PP^r}(-r)],$$
where the complex is placed in degrees $[-r,0]$, the differential is induced
by the $A_V$-action on $M$ via the natural isomorphisms
$\Hom(\Om^i_{\PP^r}(-r+i),\Om^{i-1}_{\PP^r}(-r+i-1))\simeq V^*$, and the brackets denote the convolution of the complex of objects in $D^{[0]}(X\times\PP^r)_0$.
One can easily check that $\Psi$ is a fully faithful embedding and that $\Phi\circ\Psi=\id$.
Since the image of $\Psi$ generates $D(X\times\PP^r)$ as a triangulated category,
it follows that $\Psi$ is essentially surjective, hence it is an equivalence.
 
Let us consider the following collection of functors
$$\begin{array}{rcl}T_i:\cC_V&\ra&\cC_V\\ F&\mapsto& H^i_0(F(1)),\end{array}$$
where we view $\cC_V$ as a subcategory in $D(X\times\PP^r)$, and where
$i\le 0$ (for $i>0$ the analogous functors are zero).

Now suppose $F$ is in the heart of the $0$-th $t$-structure and $M = \Phi
(F)$. The Koszul complex $$0\ra\cO_{\PP^r}\ra \big( \bigwedge^rV\big)\otimes \cO_{\PP^r}(1)\ra\ldots\ra
\big( \bigwedge^2V\big)\otimes \cO_{\PP^r}(r-1)\ra V\otimes \cO_{\PP^r}(r)\ra\cO_{\PP^r}(r+1)\ra 0.
$$ induces a quasi-isomorphism 
$$ \big[K^{\bullet}(M)\big] \to  \bR p_*(F(r+1))$$
Here we define the complex of objects in $\cC$
$$K^{\bullet}(M) \quad = \quad M_0\ra \big(\bigwedge^rV\big)\otimes
M_1\ra\ldots\ra\big(\bigwedge^2V\big)\otimes M_{r-1}\ra V\otimes M_r,$$ positioned in
degrees $[-r,\ldots,0]$, and $\big[K^{\bullet}(M)\big]$ is its
convolution. In
particular the cohomology objects of $\bR 
p_*(F(r+1))$ in terms of the $t$-structure on $X$ are the cohomology
objects of the complex $K^{\bullet}(M)$.

Write $$M_{r+1}^i = H^i(K^{\bullet}(M)),$$ in particular we write explicitly 
$$M_{r+1}^0\ \  =\ \  \Coker\Big(\big(\bigwedge^2V\big)\otimes
M_{r-1}\ \ra\ V\otimes M_r\Big).$$ Note   
that this involves only $M_{r-1}$ and $M_r$.
This discussion implies the first part of the following lemma:

\begin{lemma}\label{quiverlem} 
\begin{enumerate}\item[(i)]  For an object $M=(M_0,\ldots,M_r)$ of $\cC_V$
one has 
$$\begin{array}{rcll}T_0(M) &\simeq&
(M_1,\ldots,M_r,M_{r+1}^0),&\\ 
T_i(M)& \simeq& (\,\ 0\ ,\ldots,\,\ 0\ ,M_{r+1}^i),\ & i<0.\end{array}$$
\item[(ii)] For $M=(M_0,\ldots,M_r)\in \cC_V$ and $n>0$ one has
$$\Phi H^0_0(\Phi^{-1}M(n+r))=(M'_0,\ldots,M'_r),$$ where
$$M'_i=\Coker\Big(S^{n+i-1}(V)\otimes\big(\bigwedge^2V\big)\otimes M_{r-1}\ \
\ra\ \  
S^{n+i}(V)\otimes M_r\Big)$$
\end{enumerate}
\end{lemma}

{\bf Proof.} We need to prove  (ii). This follows from the exact sequence
(\ref{prexactseq}) for $m=0$ together 
with an observation that in this sequence one has $V_r=S^n(V)$, 
$$V_{r-1}=\ker\Big(S^n(V)\otimes V\ \ra\ S^{n+1}(V)\Big)\ \ =\ \
\im\Big(S^{n-1}(V)\otimes\big(\bigwedge^2V\big)\ \ra\ S^n(V)\otimes V\Big).$$
\qed

\begin{remark} It is easy to see that for every $n\ge 1$ and every
$F\in\cC_V$ one has 
$$H^0_0(F(n))\simeq T_0^n(F),$$
where $T_0^n=T_0\circ\ldots\circ T_0$ is the $n$-th iteration of
$T_0$. However, we will 
not use this fact.
\end{remark}

{\bf Proof of Proposition \ref{stabprop}.}

{\sc Step 1: reduction to the case $i=0$ and $F\in D^{\le 0}(X\times\PP^r)_0$.} 
According to Lemma \ref{inclem} for
$m<n$ we have $D^{[a,b]}(X\times\PP^r)_m\sub D^{[a-r,b]}(X\times\PP^r)_n$. Hence,
for every $F\in D(X\times\PP^r)$ there exists an interval $[a,b]$
such that $H^i_n(F)=0$ for $i\not\in[a,b]$ and $n\geq 0$. Let us show using
induction in $b-a$ that
we can reduce ourselves to the particular case $i=b$ (and hence to the case $i=b=0$).
The case $b-a=0$ is clear. If $b-a>0$ then by the assumption we can find $N>0$ such that
for $n\ge N$ the morphisms $H_n^b F\to H_{n+1}^b F$ are isomorphisms. Replacing
$F$ by $F(N)$ we can assume that these maps are isomorphisms for $n\ge 0$.
This implies that the natural morphisms $\tau_{<b}^nF\to\tau_{<b}^{n+1}F$ are
also isomorphisms for $n\ge 0$. Hence, we have $H^i_n(\tau_{<b}^0F)\simeq 
H^i_n(\tau_{<b}^nF)=0$ for $n\ge 0$, $i\not\in[a,b-1]$. It remains to apply
the induction assumption to $\tau_{<b}F$. 

{\sc Step 2: restatement in terms of $0$-th $t$-structure.} The inclusion
$D^{\le -1}(X\times\PP^r)_n\sub  
D^{\le -1}(X\times\PP^r)_{n+1}$ implies that the  
morphism $H^0_n(F)\ra H^0_{n+1}(F)$ can be identified with the canonical
morphism 
$H^0_n(F)\ra \tau_{\ge 0}^{n+1}(H^0_n(F))$, 
Therefore,
it is enough to show that $H^0_n(F)$ lies in the heart of the $(n+1)$-st
$t$-structure 
for all sufficiently large $n$. Using the isomorphism $H^0_n(F)(n)\simeq
H^0_0(F(n))$ 
we can restate this as the assertion that $H^0_0(F(n))(1)$ lies in the heart of
the 
$0$-th $t$-structure for $n\gg 0$. 

{\sc Step 3: restatement in terms of $M_{r-1}$ and $M_r$.} In terms of  the 
functors $(T_i)$ introduced above 
this means that $T_i(\Phi H^0_0(F(n)))=0$ for $i<0$ and $n\gg 0$.
Let $\Phi H^0_0(F)=(M_0,\ldots,M_r)\in\cC_V$. Then by Lemma \ref{quiverlem}(ii) we
have 
$\Phi H^0_0(F(n+r))=M^n=(M^n_0,\ldots,M^n_r)$, where
$$M^n_i=\Coker\Big(S^{n+i-1}(V)\otimes\big(\bigwedge^2V\big)\otimes M_{r-1}\ \
\ra\ \ 
S^{n+i}(V)\otimes M_r\Big).$$
By Lemma \ref{quiverlem}(i) we have to check that $H^{<0}(K^{\bullet}(M^n))=0$
for 
$n\gg 0$. 

{\sc Step 4: reduction in terms of modules of  $S(V)$.} 
Let us consider the natural truncation functor $M\ra M_{[0,r]}$ from the
category 
of graded $S(V)$-modules in $\cC$ to $\cC(V)$ that leaves only graded
components of $M$ in the 
range $[0,r]$. For a graded $S(V)$-module $M$ in $\cC$ set
$K^{\bullet}(M)=K^{\bullet}(M_{[0,r]})$. 
Then we have $K^{\bullet}(M^n)=K^{\bullet}(C(n))$, where $C$ is the following
graded $S(V)$-module in $\cC$:
$$C=\Coker(S(V)\otimes\big(\bigwedge^2V\big)\otimes M_{r-1}(-1)\ra
S(V)\otimes M_r).$$
Now for every graded $S(V)$-module $N$ in $\cC$ let us denote by
$\wt{K}^{\bullet}(N)$ the  
following complex (placed in degrees $[0,r+1]$):
$$0\ra N_0\ra \big(\bigwedge^rV\big)\otimes N_1\ra\ldots\ra\big(\bigwedge^2V\big)\otimes
N_{r-1}\ra V\otimes N_r\ra 
N_{r+1}\ra 0.
$$
It is enough to prove that $\wt{K}^{\bullet}(C(n))$ is exact for $n\gg 0$.

 We claim that actually for every graded $S(V)$-module $N$ of finite type in
$\cC$ 
the complex $\wt{K}^{\bullet}(N(n))$ is exact for $n\gg 0$.
 This is clear for {\em free} modules of finite type by the exactness of the
 Koszul 
 complex 
(in all degrees but one). 

{\sc Step 5: proof of exactness.}
We now prove 
by descending induction in $i$,  that for a module $N$ of finite type
one has 
$H^i(\wt{K}^{\bullet}(N(n))$ for $n\gg 0$. For $i>r+1$ this is trivial.
Assume that the statement is true for $i+1$ and let us consider an exact
sequence in the category of graded $S(V)$-modules
$$0\ra N'\ra P\ra N\ra 0,$$ 
where $P$ is a free module of finite type.
 Then it induces an exact sequence of complexes
$$0\ra\wt{K}^{\bullet}(N'(n))\ra\wt{K}^{\bullet}(P(n))\ra
\wt{K}^{\bullet}(N(n))\ra 0$$ 
and hence a long exact sequence of cohomology
$$\ldots\ra H^i(\wt{K}^{\bullet}(P(n))\ra H^i(\wt{K}^{\bullet}(N(n))\ra
H^{i+1}(\wt{K}^{\bullet}(N'(n))\ra\ldots$$
For $n\gg 0$ we have $H^i(\wt{K}^{\bullet}(P(n))=0$. Also, by part (i) of
Theorem \ref{catHilbthm} 
the module $N'$ is of finite type. Hence, by the induction assumption   
we have $H^{i+1}(\wt{K}^{\bullet}(N'(n))=0$ for $n\gg 0$ and from the above
exact sequence 
we get that $H^i(\wt{K}^{\bullet}(N(n))=0$ for $n\gg 0$.
\qed

\paragraph{\bf C. The constant sheaf of $t$-structures on $\PP^r$.}

\begin{theorem}\label{projthm}
Assume that the heart $\cC = D(X)^{\leq 0}\cap D(X)^{\geq 0}$ is
  noetherian.

  Then 
\begin{enumerate} 
\item the  subcategories
$$D(X\times\PP^r)^{\le 0}=\{F:\ \bR p_*(F(n))\in D^{\le 0}
\text{ for all } n \gg 0\},$$
$$D(X\times\PP^r)^{\ge 0}=\{F:\ \bR p_*(F(n))\in D^{\ge 0}
\text{ for all } n \gg 0\},$$
define a nondegenerate $t$-structure on $D(X\times\PP^r)$,
which extends to a sheaf of $t$-structures over $\PP^r$.
Moreover, 
\item 
the heart of this $t$-structure on $D(X\times\PP^r)$ and on  each
$D(X\times U)$ is noetherian.
\end{enumerate}
\end{theorem}

{\bf Proof.} 

\noindent{\sc Proof of (1).} Lemma \ref{inclem} easily implies that
$$D^{\le 0}(X\times\PP^r)=\cup_n D^{\le 0}(X\times\PP^r)_n,$$
$$D^{\ge 0}(X\times\PP^r)=\cap_n D^{\ge 0}(X\times\PP^r)_n.$$
This implies that $D^{\ge 1}(X\times\PP^r)$ is the right orthogonal
of $D^{\le 0}(X\times\PP^r)$. Let $F\in D(X\times\PP^r)$
be any object. By Proposition \ref{stabprop} there exists
$n$ such that all objects $(H^i_n(F), i\in\ZZ)$ belong to the hearts of all
the $t$-structures corresponding to $m>n$. It follows that
$\tau^n_{\le 0}(F)\in D^{\le 0}(X\times\PP^r)$ and
$\tau^n_{\ge 1}(F)\in D^{\ge 1}(X\times\PP^r)$ for $n\gg 0$. This allows us to
define a right adjoint $\tau_{\le 0} : D(X\times \PP^r) \to  D^{\le
0}(X\times\PP^r)$ by $ \tau_{\le 0}(F ) = \tau^n_{\le 0}(F)$ for $n\gg 0$. This
shows that this is a $t$-structure.

By the definition, this $t$-structure on $D(X\times\PP^r)$ is preserved
by the functor of tensoring with $\cO_{\PP^r}(1)$. Therefore, by Theorem
\ref{sheaflinethm} we obtain a sheaf of $t$-structures over $\PP^r$. This
proves (1).

\noindent{\sc Proof of (2) for $\PP^r$.}

Let us denote by $\cC$ the heart of the original $t$-structure on $D(X)$
and by $\cC_{\PP^r}$ the heart of the constructed $t$-structure on
$D(X\times\PP^r)$. 
With every $F\in\cC_{\PP^r}$ we can associate a graded $S(V)$-module $M(F)$ in
$\cC$ 
defined by $M(F)_n=H^0\bR p_*(F(n))$ for $n>0$ and $M(F)_n=0$ for $n\le 0$,
where  
$H^0:D(X)\ra\cC$ is the cohomology functor associated with the $t$-structure
on $D(X)$. 

\begin{claim}
 For every $F\in\cC_{\PP^r}$ the $S(V)$-module $M(F)$ is of finite
type.
\end{claim} 
{\bf Proof of claim.} 
By the definition we have $\bR p_*(F(n))\in\cC$ for $n\gg 0$, so for
such $n$ we 
have $M(F)_n=\bR  p_*(F(n))$. It suffices to show that $V \otimes M(F)_n \to
M(F)_{n+1}$ is surjective. Thus, the exact sequence 
$$0\ra\Om^1_{\PP^r}(n)\ra V\otimes\cO_{\PP^r}(n-1)\ra\cO_{\PP^r}(n)\ra 0$$
shows that it is enough to check that $\bR  p_*(F\otimes\Om^1_{\PP^r}(n))\in
D^{\le 0}$. 

Assume now that $n$ is large enough so that $\bR  p_*(F(m))\in\cC$ for $m\ge
n-r-1$. 
Using the exact sequences
$$0\ra\Om^i_{\PP^r}(n)\ra\bigwedge^i{V}\otimes\cO_{\PP^r}(n-i)\ra
\Om^{i-1}_{\PP^r}(n)\ra  0$$ 
for $i=2,\ldots,r$ we see that it is sufficient to show that
$\bR  p_*(F\otimes\Om^r_{\PP^r}(n))=\bR  p_*(F(n-r-1))$ is in $D^{\le 0}$. But
this follows from 
the assumption that $\bR  p_*(F(n-r-1))$ is in $\cC$. So the claim is proven.

We proceed to prove that the category $\cC_{\PP^r}$ is noetherian.
First note that for any exact sequence 
$$0\ra F_1\ra F_2\ra F_3\ra 0$$
in $\cC_{\PP^r}$ the induced sequence 
$$0\ra M(F_1)_n\ra M(F_2)_n\ra M(F_3)_n\ra 0$$
in $\cC$ is exact for $n\gg 0$. Next, we claim that if for some
$F\in\cC_{\PP^r}$ we have 
$M(F)_n=0$ for $n\gg 0$ then $F=0$. Indeed, this would mean that $\bR
p_*(F(n))=0$ for $n\gg 0$ 
which clearly implies $F=0$ by considering cohomology sheaves with respect to
the standard $t$-structure. 
Therefore, if for some subobject $F'\sub F$ of
an object in $\cC_{\PP^r}$ we have $M(F')_n=M(F)_n$ for $n\gg
0$ then $F'=F$.

Recall that by Theorem \ref{catHilbthm} the category of graded
$S(V)$-modules of finite type in $\cC$ is noetherian. Since as we have shown
above 
for every $F\in\cC_{\PP^r}$ the $S(V)$-module $M(F)$ is of finite type, we
immediately 
conclude that every object of $\cC_{\PP^r}$ is noetherian. 

\noindent{\sc Proof of (2) for open $U \subset \PP^r$.}

We want to prove that the heart $\cC_U$ of the $t$-structure on $X\times U$ is
noetherian. Fix an increasing sequence $F_1\subset F_2 \cdots \subset F$ of
objects in $\cC_U$. By the sheaf axiom it suffices to show that the sequence
stabilizes on a covering family of open subsets of $U$, so we may assume $U$
affine. The divisor  $\PP^r \setmin U = D$ is defined by a section $f$ of
$\cO_{\PP^r}(d)$ for some $d$. 

By \cite{TT} there is an extension $\bar F\in D(X\times \PP^r)$ of
$F$. Note that the restriction  $(H^0\bar F)_U = F$ by the
definition of a sheaf of $t$-structures, so replacing $\bar F$ by
$H^0\bar F$ we may assume $\bar F \in \cC_{\PP^r}$. We similarly
construct $\bar F_i$. Replacing $\bar F_i$ by $\bar F_i(-k_i D)$ we
may assume that the injection $F_i \hra F$ extends to a morphism $\phi_i:\bar
F_i \to \bar F$ (see Lemma \ref{Lem:extend-morphism}). Replacing
$\bar F_i$ by $$\sum_{j=1}^i \phi_j(\bar F_j)\subset 
\bar F$$ we may assume $$\bar F_1 \subset \bar F_2 \subset \cdots
\subset \bar F$$ is increasing.  This sequence stabilizes, therefore
the original sequence stabilizes, as required. 
\qed

The above proof  also leads to  the following.

\begin{lemma}\label{surjprojlem} In the above notations for every object
$F\in\cC_{\PP^r}$ there exists a surjection of the form $p^*G(n)\ra F$ in
$\cC_{\PP^r}$, 
where $G\in\cC$, $n\in\ZZ$.
\end{lemma}

{\bf Proof.} Without loss of generality we can assume that $\bR  p_*(F(n))\in\cC$
for $n>0$. 
As we have seen in the proof of the above theorem, 
the graded $S(V)$-module $M(F)$ is of finite type. Hence, there exists $N>0$
such that 
the natural map 
$$\bigoplus_{n=1}^N\ \ S(V)\ \otimes\ M(F)_n(-n)\quad \lrar\quad  M(F)$$
of graded $S(V)$-modules in $\cC$ is surjective. This easily implies that the
canonical 
map 
$$\bigoplus_{n=1}^N\ \ \Big(p^*\,\bR  p_*(F(n))\Big)(-n)\quad \lrar\quad  F$$
is surjective in $\cC_{\PP^r}$. Therefore, we have a surjection
$$\bigoplus_{n=1}^N \ \ S^{N-n}(V) \ \otimes \
\Big(p^*\,\bR  p_*(F(n))\Big)(-N)\quad \lrar\quad F.$$ 
\qed

\begin{lemma}\label{pushforlem}
The functor $\bR  p_*:D(X\times\PP^r)\ra D(X)$
is left $t$-exact.
\end{lemma}

{\bf Proof.} This follows from the fact that $\bR  p_*$ is right adjoint to
the $t$-exact functor $p^*$. \qed

\paragraph{\bf D. The constant sheaf of $t$-structures on $\PP^r\times \PP^s$.}

\begin{proposition}\label{projprodprop}
Consider the $t$-structure on $D(X\times\PP^r\times\PP^s)$
obtained by applying the construction of Theorem \ref{projthm} twice:
first
to get a $t$-structure on $D(X\times\PP^r)$ and then to get a $t$-structure
on $D(X\times\PP^r\times\PP^s)$. Then this $t$-structure can be described as
follows: $F\in D^{\le 0}$ (resp., $F\in D^{\ge 0}$) 
if and only if there exists an integer
$N$ such that
for all $m,n>N$ one has $\bR  p_*(F(m,n))\in D^{\le 0}$ (resp., $\bR  p_*(F(m,n))\in D^{\ge 0}$).
\end{proposition}

{\bf Proof.} It suffices to check that for every object
$F$ in the heart of our $t$-structure on the product $D(X\times\PP^r\times\PP^s)$
there exists $N$ such that for all $m,n>N$ one has 
$$\bR  p_*(F(m,n))\quad \in\quad \cC\ =\ D(X)^{\le 0}\cap D(X)^{\ge 0}.$$
Let $p_1:\ X\times\PP^r\times\PP^s\ \ra\ X\times\PP^r$ be the projection.
By the definition, there exists an integer $N'$ such that for all $n\ge N'$ one has
$\bR  p_{1\,*}(F(0,n))\in\cC_{\PP^r}$. By Lemma \ref{pushforlem} this implies
that $\bR  p_*(F(m,n))\in D(X)^{\ge 0}$ for all $n\ge N'$ and all $m$.
On the other hand, for every $n\ge N'$ there exists an integer $M_n$
such that for all $m\ge M_n$ one has $\bR  p_*(F(m,n))\in\cC$. Set
$M=\max(M_{N'},M_{N'+1},\ldots,M_{N'+s})$. 

We claim that for every $n\ge N'$ and
every $m\ge M$ one has $\bR  p_*(F(m,n))\in \cC$. Indeed, it suffices
to check that in this range of $(m,n)$ one has $\bR  p_*(F(m,n))\in D(X)^{\le 0}$.
We can argue by induction in $n$. For $N'\le n\le N'+s$ the assertion is true by
the choice of $M$. Now let $n>N'+s$ and assume that the assertion
is true for all $n'$ such that $N'\le n'<n$. Then the exact sequence
$$0\ra\cO_{\PP^s}(n-s-1)\ra\ldots\cO_{\PP^s}(n-1)^{s+1}\ra\cO_{\PP^s}(n)\ra 0$$
implies that $\bR  p_*(F(m,n))\in D(X)^{\le 0}$. \qed

\begin{corollary}\label{projprodcor} The two $t$-structures on
  $D(X\times\PP^r\times\PP^s)$ 
obtained from the $t$-structure on $D(X)$ by iterating the construction of
Theorem \ref{projthm} coincide. This $t$-structure extends
to a sheaf of $t$-structures over $\PP^r\times\PP^s$.
\end{corollary}

{\bf Proof.} The second assertion follows from Theorem \ref{sheaflinethm}
since our
$t$-structure on $D(X\times\PP^r\times\PP^s)$ is stable under tensoring
with the pull-back of any line bundle on $\PP^r\times\PP^s$. \qed

We also record the following lemma for later use.

\begin{lemma}\label{finsuplem}
Let $Y\subset \PP^r\times\PP^s$ be a closed subset such that
the natural projection $\pi:Y\ra \PP^r$ is a finite morphism and let
$F\in D(X\times\PP^r\times\PP^s)$ be an object supported on $X\times Y$.
Then $F\in D(X\times\PP^r\times\PP^s)^{\le 0}$ 
(resp., $F\in D(X\times\PP^r\times\PP^s)^{\ge 0}$)
if and only if $\bR  p_{1\,*}F\in D(X\times\PP^r)^{\le 0}$ 
(resp., $\bR  p_{1\,*}F\in D(X\times\PP^r)^{\ge 0}$),
where $p_1:X\times\PP^r\times\PP^s\ra X\times\PP^r$ is the projection.
\end{lemma}

{\bf Proof.} The condition that $F\in \cC_{\PP^r\times\PP^s}$ is
equivalent to $\bR  p_{1*}(F(0,n))\in \cC_{\PP^r}$ for all $n\gg 0$.
Replacing $Y$ by its finite thickening we can assume that $F$ is a push-forward
of an object in $D(X\times Y)$.
Let $(U_i)$ be an open covering of $\PP^r$ such that the line bundle 
$\cO_Y(0,1)$ is trivial on $\pi^{-1}(U_i)$. 
Then $\bR  p_{1*}(F(0,n))|_{X\times U_i}\simeq \bR  p_{1*}(F)|_{X\times U_i}$ for all $n\in\ZZ$.
It remains to use the fact that we have a sheaf of $t$-structures over
$\PP^r$.  \qed 

\subsection{The essential image of push-forward along a closed embedding}

\begin{lemma}\label{suphomlem} 
Let $i:Z\ra Y$ be the zero locus of a non-zero function $f$ on an
integral quasi-projective variety $Y$. Let also $Z'$ be the zero locus of 
$f^2$, 
$j:Z\ra Z'$ and $i':Z'\ra Y$ be the natural embeddings.
Then for every $F,G\in D(Z)$ and every morphism
$\alpha:i_*F\ra i_*G$ in $D(Y)$ there exists 
a morphism $\beta:j_*F\ra j_*G$ in $D(Z')$, such that $\alpha=i'_*\beta$.
\end{lemma}

{\bf Proof.} Let $F^{\bullet}$, $G^{\bullet}$
be (bounded) complexes of coherent $\cO_Z$-modules
representing $F$ and $G$.   We can choose 
a bounded complex $P^{\bullet}$  of torsion-free coherent 
$\cO_Y$-modules, a surjective morphism of complexes of $\cO_Y$-modules 
$q:P^{\bullet}\ra F^{\bullet}$ which is a quasi-isomorphism, and a chain map
$p:P^{\bullet}\ra G^{\bullet}$ such that
 $\alpha=pq^{-1}$.
 Hence, if $K^{\bullet}=\Ker(q)$ then $K^{\bullet}$ 
is an acyclic complex and the multiplication by $f$ is injective
on each $K^i$. We have $fP^{\bullet}\sub K^{\bullet}\sub P^{\bullet}$.
Set $\ov{P}^{\bullet}=P^{\bullet}/fK^{\bullet}$. Then both maps $q$ and $p$
factor through maps $\ov{q}:\ov{P}^{\bullet}\ra F^{\bullet}$ and
$\ov{p}:\ov{P}^{\bullet}\ra G^{\bullet}$. 
Moreover, $\ov{q}$ is still a quasi-isomorphism,
so $\alpha=\ov{p}\ov{q}^{-1}$. It remains to note that $\ov{P}^{\bullet}$
is a complex of $\cO_{Z'}$-modules. \qed

\begin{theorem}\label{divisorthm}
Let $Y_0$ be a smooth affine variety, $f$ be a function on
$Y_0$ such that its divisor of zeros $Z_0=Z(f)$ is smooth. Denote $Y=
X \times Y_0$ and $Z = X \times Z_0$. 
Denote by $i:Z \ra Y$ the natural embedding.
Let $F\in D(Y)$,  and assume that the morphism $F\ra F$
given by the multiplication with $f$ is zero. Then
$F\simeq i_*F'$ for some $F'\in D(Z)$.
\end{theorem}

{\bf Proof.} By the projection formula we have
$i_*\bL i^*F\simeq i_*\cO_{Z}\otimes F$. 
Now the exact sequence 
$$0\ra \cO_Y\stackrel{f}{\ra}\cO_Y \ra i_*\cO_{Z}\ra 0$$
gives rise to an exact triangle
$$F\stackrel{f}{\ra} F\ra i_*\cO_Z\otimes F\ra F[1].$$
Therefore, if the map $F\stackrel{f}{\ra} F$ is zero then
$F$ is a direct summand in $i_*\bL i^*F$. Let 
us denote by $\pi:i_*\bL i^*F\ra F$ the corresponding morphism, 
so that $\pi\circ\can_F=\id_F$, where $\can_F:F\ra i_*\bL i^*F$
is the canonical morphism.

 Since the cohomology sheaves $H^iF$ (with respect to the standard
 $t$-structure) are push-forwards of
some sheaves on $Z$ there exists a finite order thickening $Z'$ of $Z$
in $Y$ such that $F\simeq i'_*G$ for some $G\in D(Z')$,
where $i':Z'\ra Y$ is the natural embedding. Let us denote also
by $j:Z\ra Z'$ the natural embedding, so that $i=i'\circ j$.
Then we have $i_*\bL i^*F\simeq i'_*j_*\bL i^*F$, so $\pi$ (resp., $\can_F$) 
is an element of $\Hom(i'_*j_*\bL i^*F, i'_*G)$
(resp., $\Hom(i'_*G,i'_*j_*\bL i^*F)$).

 By Lemma \ref{suphomlem},
replacing $Z'$ by a thickening we may assume that
$\pi=i'_*\phi$, $\can_F=i'_*\psi$ for some morphisms $\phi:j_*\bL i^*F\ra G$
and $\psi:G\ra j_*\bL i^*F$ in $D(Z')$. 
Moreover, since $i'_*(\phi\circ\psi)=\id_F$,
it follows that $\phi\circ\psi:G\ra G$ induces an isomorphism on
cohomology sheaves, hence it is an isomorphism.

Since $Z_0$ is smooth and affine there exists a morphism
$p:Z'\ra Z$ such that $p\circ j=\id$. Then we can
consider the morphism $p_*\phi:\bL i^*F\simeq p_*j_*\bL i^*F\ra p_*G$
in $D(Z)$.
Finally, we define the morphism $\alpha:G\ra j_*p_*G$ in $D(Z')$ 
by setting $\alpha\, :=\, (j_*p_*\phi)\circ\psi$:
$$\xymatrix{G\ar[r]^{\psi}\ar[rd]_{\alpha} & j_*\bL i^*F\ar[d]^{j_*p_*\phi}\\&
j_*p_*G.}$$ 
We claim that $\alpha$ is an isomorphism. Since the morphism $p$
is finite, it suffices to check that $p_*\alpha$ is an isomorphism.
But $p_*\alpha=p_*\phi\circ p_*\psi=p_*(\phi\circ\psi)$ and
$\phi\circ\psi$ is an isomorphism. Setting $F' = p_*G$ we obtain 
$$ i_*F' \simeq i_*p_*G \simeq i'_*j_*p_*G \simeq i'_* G \simeq F, $$
 as required. \qed 

\subsection{Restricting a sheaf of $t$-structures to a closed subset}

\begin{lemma}\label{isommorlem} 
Let $i:Z\ra Y$ be a regular embedding of codimension $r$,
where $Z=Z(f_1,\ldots,f_r)$ is the zero locus of an $r$-tuple of functions
$f_1,\ldots,f_r$ on $Y$. Let $F,\, F'\ \in\ D(Z)$ be two objects such that
$\Hom^n(i_*F,i_*F')=0$ whenever $n<0$. Then 
$\Hom^n(F,F')=0$ whenever $n<0$, and the natural map
$\Hom(F,F')\ra\Hom(i_*F,i_*F')$ is an isomorphism.
\end{lemma}

{\bf Proof.} It suffices to consider the case when $Z$ is a divisor in $Y$:
$Z=Z(f)$ for a function $f$ on $Y$.
Note that $\Hom^n(F,F')=0$ for $n \ll 0$, so we can prove the
first assertion by induction. 
Assume that we know that $\Hom^{n'}(F,F')=0$ for
$n'<n$, where $n<0$.
The canonical distinguished triangle
$$F[1]\ra \bL i^*i_*F\ra F\ra F[2]$$
induces a long exact sequence
$$\ldots\ra\Hom^{n-2}(F,F')\ra\Hom^n(F,F')\ra\Hom^n(i_*F,i_*F')\ra
\Hom^{n-1}(F,F')\ra\ldots$$
Since $\Hom^{n-2}(F,F')=\Hom^{n-1}(F,F')=0$ by the induction assumption,
we derive that $\Hom^n(F,F')\simeq\Hom^n(i_*F,i_*F')=0$ which
finishes the proof of the first assertion.
Now using the same exact sequence for $n=0$ we derive the second
assertion. \qed

\begin{theorem}\label{restrictthm} 
Let $S$ be a smooth variety, $T\sub S$ be a smooth closed
subvariety. Then for every sheaf of $t$-structures over $S$ there exists
a unique sheaf of $t$-structures over $T$ such that the push-forward
functor $D(X\times T)\ra D(X\times S)$ is $t$-exact.
\end{theorem}

{\bf Proof.} Let $i:X \times T\ra X \times S$ denote the embedding. We claim that
$$D(X\times T)^{\le 0}=\{F:\ i_*F\in D^{\le 0}\},$$
$$D(X\times T)^{\ge 0}=\{F:\ i_*F\in D^{\ge 0}\},$$
is a nondegenerate $t$-structure on $D(X\times T)$.
Let us first check that for $A\in D(X\times T)^{\le 0}$,
$B\in D(X\times T)^{\ge 1}$ one has $\Hom(A,B)=0$.
It suffices to prove that $\bR  \cHom_T(A,B)\in D(T)^{\ge 1}$
(we use the notation of \ref{sheafsec}).
This assertion is local on $T$, so we can assume that $S$ is affine
and that $T$ is the zero locus of an $r$-tuple of functions $(f_1,\ldots,f_r)$
on $S$ (where $r=\dim S-\dim T$). Since $T$ is affine it suffices to
check that $\Hom^i(A,B)=0$ for $i\leq 0$ in this case. But this 
follows immediately from Lemma \ref{isommorlem}.

It remains to check
that $D(X\times T)$ is generated as a triangulated category by
$\cC_T:=D(X\times T)^{\le 0}\cap D(X\times T)^{\ge 0}$.
Let $F\in D(X\times T)$ be an arbitrary object. Without loss
of generality we can assume that $i_*F\in D(X\times S)^{\ge 0}$ but
$i_*F\not\in D(X\times S)^{\ge 1}$. Let us consider the canonical
morphism $\phi:\tau_{\le 0}i_*F\ra i_*F$. It suffices to prove that
there exists an object $F^0\in D(X\times T)$ and a morphism
$\ov{\phi}:F^0\ra F$ such that
\begin{equation}\label{truncisom}
\tau_{\le 0}i_*F\simeq i_*F^0
\end{equation}
and $\phi=i_*\ov{\phi}$.
Indeed, then we would consider the cone $F'$ 
of the morphism $\ov{\phi}$. Since $i_*F'$ has smaller cohomological
length than $i_*F$ with respect to our $t$-structure, we would be
able to finish the proof using induction.

Note that since $F'\in D(X\times T)^{\ge 1}$, using orthogonality of
$D(X\times T)^{\le 0}$ and $D(X\times T)^{\ge 1}$ one can immediately
check that the pair $(F^0,\ov{\phi})$ (if exists) is uniquely determined
up to a unique isomorphism.
To prove existence let us first consider a similar problem for the
restrictions $F_i=F|_{X\times (T\cap U_i)}$,
where $(U_i)$ is sufficiently fine open affine covering of $S$.
Indeed, locally we can assume that $T$ is given in $S$ by $r$ equations
$f_1,\ldots,f_r$ (such that all intersections of divisors $(f_i)$ are
smooth). Furthermore,
we can reduce to the case when $T=Z(f)$ is a divisor in
$S$. In this case using the universal property of $\tau_{\le 0}$ one
immediately checks that multiplication by $f$ gives a
zero endomorphism of $\tau_{\le 0}i_*F$, so
we can apply Theorem \ref{divisorthm} to find $F^0\in D(X\times T)$
such that (\ref{truncisom}) holds.
By Lemma \ref{isommorlem} there exists a unique morphism
$\ov{\phi}:F^0\ra F$ such that $\phi=i_*\ov{\phi}$.
Applying this argument to all restrictions
$F_i\in D(X\times (T\cap U_i))$
we find a collection of objects $F^0_i\in D(X\times (T\cap U_i))$
and morphisms $\ov{\phi}_i:F^0_i\ra F_i$ with required property.
By uniqueness, we get the gluing data for these pairs. 
It remains to note that since $\Hom^j(F^0_i,F^0_i)=0$ for $j<0$
we have $\bR  \cHom_{U_i}(F^0_i,F^0_i)\in D(U_i)^{\ge 0}$. Therefore,
we can apply Corollary \ref{sheafcor} to glue objects $(F^0_i)$
into an object $F^0\in D(X\times T)$. Similarly,
we have $\bR  \cHom_{U_i}(F^0_i,F_i)\in D(U_i)^{\ge 0}$. Hence, by
Lemma \ref{sheaflem} we can glue morphisms $\ov{\phi}_i$ into
a morphism $\ov{\phi}:F^0\ra F$. \qed

\begin{lemma}\label{closedemblem} 
In the situation of the above theorem the functor $\bL i^*:D(X\times
S)\ra D(X\times T)$
is right $t$-exact and for every $F\in \cC_S :=  D^{\le 0}(X\times
S)\cap D^{\ge 0}(X\times S)$ 
one has 
$H^0\bL i^*i_*F\simeq F$.
\end{lemma}

{\bf Proof.} The functor $\bL i^*$ is right $t$-exact since it is left adjoint to the
$t$-exact functor $i_*$, so it remains to check the second assertion of the lemma.
Since we are dealing with sheaves of $t$-structures, the assertion is local.
Also, if $i$ decomposes into a composition $i_1\circ i_2$ of closed embeddings
(and the intermediate variety is smooth) then it suffices to prove the lemma
for $i_1$ and $i_2$. Therefore, it is enough to consider the case when $T$ is a divisor
in $S$ defined as the zero locus of some function $f\in\cO(S)$. But in this case the
assertion follows immediately from the exact triangle
$$F[1]\ra \bL i^*i_*F\ra F\ra F[2].$$
\qed

\subsection{The constant sheaf of $t$-structures over a smooth projective
variety}

\begin{theorem}\label{projtstrthm} 
Assume that the heart $D(X)^{\leq 0}\cap D(X)^{\geq 0}$ is
  noetherian.

 Let $S$ be a smooth projective variety. 
Let $L$ be an ample line bundle on $S$. Then the subcategories
$$D(X\times S)^{\le 0}_L=\{F:\ \bR  p_*(F\otimes L^n)\in D^{\le 0}
\text{ for all } n \gg 0\},$$
$$D(X\times S)^{\ge 0}_L=\{F:\ \bR  p_*(F\otimes L^n)\in D^{\ge 0}
\text{ for all } n \gg 0\},$$
define a nondegenerate $t$-structure on $D(X\times S)$ with a
noetherian heart.
Furthermore, this $t$-structure does not depend on a choice of $L$
and extends uniquely to a sheaf of $t$-structures over $S$.
\end{theorem}

{\bf Proof.} Theorems \ref{projthm} and \ref{sheaflinethm}
give the required sheaf of $t$-structures over every projective space.
Now let $L$ be an ample line bundle on $S$,
$i:S\ra\PP^r$ be the closed embedding given by some power $L^n$ of $L$. 
Applying Theorem
\ref{restrictthm} we get a sheaf of $t$-structures over $S$,
extending the $t$-structure $(D(X\times S)^{\le 0}_{L^n}, 
D(X\times S)^{\ge 0}_{L^n})$ on $D(X\times S)$. 
Moreover, since the functor $i_*$ is $t$-exact and $i_*(F)=0$ only for $F=0$,
we immediately obtain that the heart of the $t$-structure on $D(X\times S)$
is Noetherian.
By Theorem \ref{sheaflinethm} the functor of tensoring with $L$ is $t$-exact.
Therefore, for every object $F$ in $D(X\times S)^{\le 0}_{L^n}$
we have $F\otimes L,\ldots, F\otimes L^{n-1}\in D(X\times S)^{\le 0}_{L^n}$,
hence, $F\in D(X\times S)^{\le 0}_L$. Thus, 
$D(X\times S)^{\le 0}_L=D(X\times S)^{\le 0}_{L^n}$. Similarly,
$D(X\times S)^{\ge 0}_L=D(X\times S)^{\ge 0}_{L^n}$. It remains to show that
this $t$-structure does not depend on a choice of projective embedding.
Assume that we have two projective embeddings $i_1:S\ra\PP^r$ and
$i_2:S\ra\PP^s$, corresponding to very ample line bundles $L_1$ and $L_2$ on $S$.
Let $i:S\ra\PP^r\times\PP^s$ be the diagonal embedding. 
Recall that by Corollary \ref{projprodcor} the two natural $t$-structures on 
$D(X\times\PP^r\times\PP^s)$ coincide. Now applying Lemma \ref{finsuplem}
we see that $F\in D(X\times S)^{\le 0}_{L_1}$ (resp., $F\in D(X\times S)^{\ge 0}_{L_1}$)
if and only if $i_*F\in D(X\times\PP^r\times\PP^s)^{\le 0}$ (resp.,
$i_*F\in D(X\times\PP^r\times\PP^s)^{\ge 0}$). Since a similar assertion holds
for $L_2$, the $t$-structures associated with $L_1$ and $L_2$ coincide. 
The obtained $t$-structure on $D(X\times S)$ is invariant with respect
to tensoring by any ample line bundle. Hence, by Theorem \ref{sheaflinethm}
it extends to a sheaf of $t$-structures. \qed

Let $\cC$ be the heart of the $t$-structure on $D(X)$ and $\cC_S$ be the heart of the
$t$-structure on $D(X\times S)$ constructed above.

\begin{lemma}\label{resolutionlem} In the situation of Theorem \ref{projtstrthm} 
every object $F$ of $\cC_S$ has a left resolution in
$\cC_S$ of the form 
$$\ldots\ra p^*G_2\otimes L^{n_2}\ra p^*G_1\otimes L^{n_1}\ra p^*G_0\otimes L^{n_0}\ra F\ra 0$$
for some integers $n_0>n_1>n_2>\ldots$ and some objects $G_i\in\cC$.
\end{lemma}

{\it Proof.} It suffices to check that for every $F$ there exists a surjection
of the form $p^*G\otimes L^n\ra F$, where $G\in\cC$, $n\in\ZZ$. 
In the case when $S$ is a projective space this follows from Lemma \ref{surjprojlem}.
In the general case we can assume that $L$ is very ample and consider the projective embedding 
$i:S\ra\PP^r$ given by $L$. There exists a surjection in $\cC_{\PP^r}$
of the form $\pi^*G(n)\ra i_*F$, where $G\in\cC$, $n\in\ZZ$, $\pi:X\times\PP^r\ra X$
is the natural projection. Since the functor 
$\bL i^*:D(X\times\PP^r)\ra D(X\times S)$ is right $t$-exact 
the morphism
$$H^0\big(\bL i^*(\pi^*G(n))\big)\ra H^0\big(\bL i^*i_*F\big)$$
is surjective in $\cC_{S}$. But $H^0\bL i^*i_*F\simeq F$ by Lemma
\ref{closedemblem} 
and $$H^0\big(\bL i^*(\pi^*G(n))\big)\ \ \simeq\ \  H^0p^*G\otimes L^n\ \
\simeq\ \ 
p^*G\otimes L^n.$$ 
\qed

\subsection{The constant sheaf of $t$-structures over a smooth
quasi-projective variety} 

In this section we use freely the fact that the category $D(V)$ is
small for any variety $V$. We also use the notion of a pre-aisle,
namely a suspended subcategory $\cS\subset \cT$ in a triangulated
category, where the shift in $\cS$ is induced from $\cT$ and the
triangles in $\cS$ come from $\cT$. See,
e.g., \cite{LST}. A pre-aisle gives a $t$-structure if and only if
there is a ``truncation functor'' $\tau:\cT\to \cS$ which is right
adjoint to the inclusion $\cS \subset \cT$. It is said to be {\em
  cocomplete} if it is stable under coproducts.

Let us denote by $D_{qc}(X)$ the unbounded derived category of quasicoherent
sheaves on $X$. The following lemma shows that a $t$-structure
$(D^{\le 0}, D^{\ge 0})$ 
on $D^b(X)$ extends to a $t$-structure on $D_{qc}(X)$.

\begin{lemma}\label{quasicohlem} 
Let $D^{\le 0}_{qc}(X)$ be the cocomplete pre-aisle in $D_{qc}(X)$
generated by  
$D^{\le 0}\sub D^b(X)$. Let also $D^{\ge 0}_{qc}(X)$ be the right orthogonal of
$D^{\le -1}\sub D^b(X)$ in $D_{qc}(X)$. Then $(D^{\le 0}_{qc}(X),
D^{\ge 0}_{qc}(X))$ 
is a $t$-structure on $D_{qc}$ and one has $D^{[a,b]}_{qc}(X)\cap
D^b(X)=D^{[a,b]}$. 
\end{lemma}

{\bf Proof.} It is easy to check that $D^{\ge 0}_{qc}(X)$ is exactly
the right orthogonal of 
$D^{\le -1}_{qc}(X)$ (see Lemma 3.1 of \cite{LST}). 
Furthermore, by Theorem A.1 of loc. cit.
these subcategories form a $t$-structure. The equality 
$$D^{\ge 0}_{qc}(X)\cap D^b(X)=D^{\ge 0}$$
is clear since both sides coincide with the right orthogonal to
$D^{\le -1}$ in $D^b(X)$. 
We have the inclusion $D^{\le 0}\sub D^{\le 0}_{qc}(X)\cap D^b(X)$ by 
the definition. 
The inverse inclusion follows from the fact that $D^{\le
0}_{qc}(X)\cap D^b(X)$ 
is left orthogonal to $D^{\ge 1}$.
\qed

We use this extension of a $t$-structure to quasicoherent complexes in
the following theorem. 

\begin{theorem}\label{Th:quasiproj-characterisation}
 Assume that the heart $\cC= D(X)^{\leq 0}\cap D(X)^{\geq 0}$ is
noetherian and {\em bounded with respect to the standard $t$-structure on
$D^b(X)$}.  
Let $S$ be a smooth quasiprojective variety. For an open subset $U\in S$ set
$$D^{\ge 0}(X\times U)=\{F\in D^b(X\times U) |\
\bR  p_*(F|_{X\times U'})\in D^{\ge 0}_{qc}(X) \text{ for every open } U'\sub
U\}$$  
and define $D^{\le 0}(U)$ to be the left orthogonal to $D^{\ge 1}(X\times U)$.
Then
\begin{enumerate} 
\item  these subcategories define a sheaf of nondegenerate $t$-structures over
$S$.  
\item When $S$ is projective this sheaf of $t$-structures coincides with the
one 
defined in Theorem \ref{projtstrthm}. 
\item When $S$ is quasi-projective, the
$t$-structure on $S$ coincides with the one induced from the sheaf of
$t$-structures on any projective completion of $S$.
\end{enumerate}
\end{theorem}

{\bf Proof.} Without loss of generality we can assume $S$ to be projective.
Let $L$ be an ample line bundle on $S$. Then for every open subset $U\sub S$
we have a $t$-structure $(D^{\le 0}(X\times U), D^{\ge 0}(X\times U))$ on
$D^b(X\times U)$ coming from Theorem \ref{projtstrthm}. Using Lemma \ref{quasicohlem}
we can extend this $t$-structure to the unbounded derived category of quasicoherent 
sheaves on $X\times U$. Let $P_L(U)\sub D_{qc}(X\times U)$ be the
cocomplete pre-aisle generated by
$(p^* D^{\le 0}\otimes L^n, n\in\ZZ)$. 

\begin{claim} The pre-aisle $P_L(U)$ coincides with
$D^{\le 0}_{qc}(X\times U)$ of Lemma \ref{quasicohlem}. 
\end{claim}
{\bf Proof of claim.} 

 Note that the inclusion $P_L(U)\sub D^{\le 0}_{qc}(X\times U)$
is clear. To check the inverse inclusion let us first consider the case $U=S$.

\noindent{\sc The case $U = S$.}
Since $D^{\le 0}_{qc}(X\times S)$ is the minimal pre-aisle containing
$D^{\le 0}(X\times S)$, 
it suffices to check that $D^{\le 0}(X\times S)\sub P_L(S)$. Moreover,
it is enough to check 
that $\cC_S\sub P_L(S)$. Using Lemma \ref{resolutionlem} we can find a
resolution of 
an arbitrary object $F\in\cC_S$ of the form 
$$\ldots\ra P_2\stackrel{d_2}{\ra} P_1\stackrel{d_1}{\ra}
P_0\stackrel{d_0}{\ra} F\ra 0$$ 
with $P_i\in\cC_S\cap P_L(S)$. 

Note that our assumption that $\cC$ is bounded
with respect to the standard $t$-structure on $D^b(X)$ immediately
implies (by taking the pushforward) that
$\cC_S$ is bounded with respect to the standard $t$-structure on
$D^b(X\times S)$.
In particular, there exists $N>0$ such that $\Hom^{>N}_{D(X\times
S)}(\cC_S,\cC_S)=0$. 
Set $K_n=\ker(d_n)\sub P_n$. Then we have a sequence of morphisms
in $D^b(X\times S)$
$$F\ra K_0[1]\ra K_1[2]\ra\ldots$$
Note that, by the definition of $N$, the composed map $F\ra K_n[n+1]$
is zero for $n>N$. 
Let $C_n\sub D^b(X\times S)$ be the convolution of the complex
$$P_n\ra P_{n-1}\ra\ldots\ra P_0$$
(this convolution exists since $P_i$ lie in the heart of some $t$-structure).
Then we have an exact triangle
$$K_n[n]\ra C_n\ra F\ra K_n[n+1].$$
Therefore, $F$ is a direct summand of $C_n$ for $n>N$. But $C_n$ belongs
to $P_L(S)$. Recall  that cocomplete pre-aisles are closed under taking 
direct summands - see \cite{LST}, Corollary 1.4: if the complement of $F$ is
$F^\perp$ then $$C_n^{\oplus \infty} \ \simeq\  \bigoplus_{i=1}^\infty (F
\oplus F^\perp)  \ \simeq\  
\bigoplus_{i=1}^\infty ( F^\perp\oplus F  ) \  \simeq\  F^\perp\oplus
\bigoplus_{i=2}^\infty (F \oplus F^\perp), $$ and the cone of the canonical
embedding $$C_n^{\oplus\infty} \  \simeq\  F^\perp\oplus \bigoplus_{i=2}^\infty
(F  \oplus 
F^\perp) \quad \lrar \quad C_n^{\oplus\infty} \ \simeq \ \bigoplus_{i=1}^\infty (F
\oplus 
F^\perp) $$ is isomorphic to $F$.  Hence, $F$ lies in $P_L(S)$. 

\noindent{\sc The general case.}
Now let us consider the case of general open subset $j:U\hra S$. 
Note that the set of $F\in D_{qc}(X\times S)$ such that $j^*F\in P_L(U)$ forms
a pre-aisle in $D_{qc}(X\times S)$. Hence, we have $j^*P_L(S)\sub P_L(U).$
Therefore,
$$D^{\le 0}(X\times U)=j^*D^{\le 0}(X\times S)\sub j^*P_L(S)\sub P_L(U)$$
and as above this implies the inclusion $D^{\le 0}_{qc}(X\times U)\sub P_L(U)$.
This finishes the proof of our claim.

We continue with the proof of the theorem.
The subcategory $D^{\ge 1}_{qc}(X\times U)$ is defined as the right orthogonal
to 
$D^{\le 0}_{qc}(X\times U)=P_L(U)$. Therefore, by Lemma 3.1 of \cite{LST} 
we have
$$D^{\ge 0}_{qc}(X\times U)=\{F\in D_{qc}(X\times U)|\
\Hom^{<0}(p^*G\otimes L^n,F)=0 \text{ for all } G\in D^{\le 0}, n\in\ZZ\}.$$
Using the isomorphisms
$$\Hom^i(p^*G\otimes L^n,F)\simeq\Hom^i(G,\bR  p_*(F\otimes L^{-n})),$$
we deduce that
$$D^{\ge 0}_{qc}(X\times U)=\{F\in D_{qc}(X\times U)|\ \bR
p_*(F\otimes L^n)\in D^{\ge 0}_{qc}
\text{ for all } n\in\ZZ\}.$$
Note that if $U$ is sufficiently small then the restriction of $L$ to
$U$ is trivial, 
hence for such $U$ we will get
$$D^{\ge 0}_{qc}(X\times U)=\{F\in D_{qc}(X\times U)|\ \bR  p_*F\in
D^{\ge 0}\}.$$ 
It remains to restrict this equality to $D^b(X\times U)$ and to use
the fact that we  
have a sheaf of $t$-structures to finish the proof.
\qed

\section{Families of objects of the heart and their
extensions}\label{Sec:families-in-heart}

\subsection{Torsion subobjects} 

As before, let $(D(X)^{\le 0}, D(X)^{\ge 0})$ be a nondegenerate
$t$-structure with noetherian heart $\cC$. We fix a smooth
projective variety $S$.
We denote the heart of  the lifted $t$-structure $(D(X\times S)^{\le
  0}, D(X\times S)^{\ge 0})$ by $\cC_S$. For an open subset $S'\subset
S$ we denote
the heart coming from the sheaf of $t$-structures $\cC_{S'\subset S}$,
or $\cC_{S'}$ if there is no possibility of confusion. Recall that we
have already shown that $\cC_{S'}$ is independent of the
compactification $S$ {\em if $\cC$ is bounded with respect to the
standard $t$-structure}, with an explicit description of this
heart. We will soon prove the independence of $\cC$ without the
boundedness assumption.

\begin{definition} 
Define an object $F\in \cC_S$ to be {\em $S$-torsion} if it is the push forward of an object $F' \in
D(X\times T)$ for some closed subscheme $T\subset
S$. Equivalently, for every divisor $D \subset S$ containing $T$, with defining
equation $f \in \cO_S(D)$, there is an integer $k$ such that the
morphism $f^k: F \to F\tototi \cO_S(kD)$ is zero.

An object $F \in \cC_S$ is {\em $S$-torsion-free} if it contains no
nonzero $S$-torsion subobject.

An object $E\in \cC_{S}$ is said to be {\em $t$-flat} if
  for every closed point $s\in S$ we have $\bL i_s^* E\ \in\ \cC$. 
\end{definition}

The following is immediate:
\begin{proposition}\label{Prop:max-torsion}
Let $E\in \cC_{S}$. Then
 \begin{enumerate} 
\item there is a maximal $S$-torsion
  subobject $F \subset E$; and
\item for every closed $T\subset S$ there is a maximal $S$-torsion
object supported over $T$.
\end{enumerate}
\end{proposition}

{\bf Proof.} We prove the first statement, the other being similar.
Consider an increasing sequence 
$$F_1 \subset F_{2}\subset\cdots \subset  E$$
 of
$S$-torsion subobjects.  Since $\cC_S$ is noetherian, this sequence
stabilizes. 
 \qed 

 \begin{corollary}
\begin{enumerate} 
\item  Let $E\in \cC_{S}$. Then there exists an $S$-torsion-free object $\bar
E\in 
 \cC_{S}$, an epimorphism $\phi:E \to \bar E$, and a dense  open immersion
 $j: S' \to S$ such that $j^*\phi: j^*E \to j^*\bar E$ is an
 isomorphism. 
\item If $D \subset S$ is a smooth divisor with defining
 section $f \in \cO_S(D)$, then $\bL i_D^*
 \bar E \in \cC_D$ and $$i_{D\,*} \bL i_D^*
 \bar E \ \  =\ \  Cone (\bar E (-D) \stackrel{f}{\to} \bar E )\ \ =\
\ \bar E/ f(\bar E(-D)).$$  
\item 
If $\dim S = 1$ then $\bar E$ is $t$-flat.
\end{enumerate}
\end{corollary}

{\bf Proof.} Let $F\subset E$ be the maximal $S$-torsion subobject,
supported over some closed $T\subset S$. Let $S' = S\smallsetminus T$.
Let $\bar E = E/F$. Then   $j^*\phi: j^*E \to j^*\bar E$ is an
isomorphism. If $\bar G \subset \bar E$ is $S$-torsion then its inverse
image $G \subset E$ is $S$-torsion as well, containing $F$, and thus
$G=F$ and $\bar G = 0$. Thus $\bar E$ is $S$-torsion free. 

For the second statement we have $\bar E\otimes
p_S^*\cO_S(-D) \to
\bar E$ injective. Thus the cone of this morphism is in
$\cC_S$. But this cone is isomorphic to $i_{D\,*}\,\bL i_D^*\,\bar E$. Since
$i_{D\,*}$ is $t$-exact and sends nonzero objects to nonzero objects,
we have that $\bL i_D^*\bar E \in \cC$. This also implies the last
statement, namely $\bar E$ is $t$-flat when $\dim S =1$.\qed.

\subsection{Extensions in the heart across a divisor}

First we note that the sheaf of $t$-structures is flasque:

\begin{lemma}\label{Lem:heart-extension} Let $j: S'\to S$ be a dense
  open immersion, and  $E'\in 
  \cC_{S'}$. Then there is an object $E\in \cC_S$ such that $j^*E =
  E'$.
\end{lemma}

{\bf Proof.} By \cite{TT} there is an object $F\in D(X\times S)$
  such that $j^*F =   E'$. Consider $E = H^0(F)$. Now by Definition
  \ref{Def:sheaf} we have $j^*E = H^0(j^*F) = H^0(E') = E'$. \qed

We have the following: 
\begin{proposition}\label{Prop:extension-across-divisor} Let $i_D:D\hra S$ be
the 
closed immersion of a 
smooth divisor, and denote by $$ S':= S\setmin D\ \ \stackrel{j}{\hra} \ \ S$$ the open immersion  of the complement.
Fix $E' \in \cC_{S'}$. Then 
\begin{enumerate}
\item  there exists an  object  $E\in \cC_S$, with no
$S$-torsion supported in $D$, such that $j^*E = E'$.
\item  If $E'$ is $S'$-torsion-free then $E$ is $S$-torsion free.
\item If, moreover, $\dim S = 1$ then $E$ is $t$-flat.
\end{enumerate}
\end{proposition} 

{\bf Proof.} By Lemma \ref{Lem:heart-extension} there is an extension
$E_0 \in \cC_S$. Let $F$ be the maximal $S$-torsion object supported
in $D$. Then $j^* (E_0/F) = E'$, and $E_0/F$ has no
$S$-torsion supported in $D$. The rest is immediate. \qed
 
\subsection{Families of objects in $\cC$}

\begin{definition}\label{Def:family-in-C} Let $S$ be a scheme of finite type
over the base field $k$.  
{\em A family of  objects in $\cC$  parametrized by  $S$}
 is an object $$E
\in D(X\times S)$$ such that for every closed point $s\in S$ we
have
$$\bL i_s^* E \in \cC.$$
\end{definition}

We have the following {\em open heart property:}

\begin{proposition} \label{Prop:open-heart}
Let $E \in D(X\times S)$ and let $T\subset S$ be a smooth closed
subscheme such that 
$\bL  i_T^*E \in \cC_T$. Then there is an open neighborhood  $T \subset U
\subset S$ such that 
$E_U \in \cC_U$. 
\end{proposition}

Applying the sheaf property we obtain:
\begin{corollary}\label{Cor:family-is-in-CS}
 Let $S$ be a smooth quasi-projective variety.
Let $E \in D(X\times S)$ be a family of  objects in $\cC$
parametrized by  $S$. Then $E \in \cC_S$.
\end{corollary}

We rely on the following lemma:

\begin{lemma}\label{Lem:H0-test} Let $S$ be a smooth quasi-projective scheme. 
Let $E\in \cC_S$, let $T\subset S$ be a smooth subscheme, and assume 
$H^0( \bL  i_T^*E) =0$. Then there is an open neighborhood $T\subset U\subset S$ such 
that $E_U =0$. 
\end{lemma}

{\bf Proof of lemma.} 
We apply induction on the codimension.

{\sc Step 1:} we consider the case where $T$ is a divisor in $S$.

We have an exact sequence in $\cC_S$: 
$$ 0 \to  F \to E \to \bar E \to 0$$
where $F$ is the maximal $S$-torsion subobject supported over $T$. Since
$\bL i_T^*$ is right $t$-exact, we have a surjection $H^0(\bL i_T^*E) \to
H^0(\bL i_T^*\bar E)$. It follows that $\bL i_T^* \bar E = H^0(\bL i_T^*
\bar E) = 0$. Considering cohomology sheaves of the standard
$t$-structure, we have that $\bar E$ is
supported away from $T$, and thus it vanishes on some open
neighborhood of $T$.  

So we may assume $E=F$ is supported over $T$. Let $f\in H^0(S,
\cO_S(T))$ be the defining section of $T$. We have that $E$ is annihilated by
$f^k$ for some minimal integer $k\geq 0$.  We apply
induction on $k$. We have that $i_{T\,*}\bL i_T^* E = Cone(E(-T) \stackrel{f}{\to}
E)$, and therefore 
$$  i_{T\,*}H^0(\bL i_T^* E)\ \  =\ \ H^0(i_{T\,*}\bL i_T^* E)\ \  =\ \ 
E/fE(-T).$$ 
Since $H^0(\bL i_T^* E) = 0$ we have $E = fE(-T)$ and therefore $E$ is
annihilated by $f^{k-1}$, contradicting minimality.

{\sc Step 2.} Consider the case $\dim S - \dim T>1$ and assumed the
result holds for lower codimensions. Write 
$$T \stackrel{i_{T\subset S_1}}{\hra} S_1 \stackrel{i_{S_1}}{\hra} S$$
 where $S_1\subset S$ is a smooth divisor. Then 
$$H^0 (\bL i_T^* E)\ \  =\ \  H^0\big(\, \bL i_{T\subset S_1}^*\, H^0(\bL
i_{S_1}^*E\,)\ \big).$$ 
By induction we have that there is an open neighborhood $T \subset U_1
\subset S_1$ such that  $$ \Big(H^0(\bL
i_{S_1}^*E)\,\Big)_{U_1} = 0.$$
Replacing $S$ by  by an open subset containing $U_1$ as a divisor, we
get   using Step 1 that $E_U=0$, which is what we needed.
\qed

{\bf Proof of Proposition.}  
We apply induction on the codimension as in the Lemma, and it suffices to
consider the case where $T$ is a divisor. Then $\bL  i_T^*$ is right exact
and $\bL  i_T^*[-1]$ is left exact.

Consider $$M = \max\{i| \Supp(\tau_{\ge i}E) \cap X\times T \neq
\emptyset\}.$$ 
Replacing $S$ by an open we may assume $\tau_{> M}E =0$ and thus
$\tau_{\ge M}E =H^M(E)[-M]$.
The distinguished triangle $\tau_{<M}E \to E \to \tau_{\ge M}E$
induces $\bL  i_T^*(\tau_{<M}E) \to \bL  i_T^* E \to \bL  i_T^*(\tau_{\ge
M}E)$. 

 Assume
by contradiction $M>0$. Since
$\bL  i_T^*$ is right exact we get that 
$$H^0\big(\bL  i_T^* (H^ME)\big) \ \ = \ \ 
H^M\big(\bL  i_T^*(\tau_{\ge M}E)\big)\ \ = \ \  H^M(\bL  i_T^*E) = 0$$ 
so, by Lemma \ref{Lem:H0-test}, $H^M(E)_U = 0$, which contradicts the
definition of $M$. 

This implies $M\leq 0$. We now replace $S$ by an open neighborhood of $T$ so
that  $\tau_{> 0}E =0$.

Taking the long exact sequence of cohomology of the distinguished
triangle $\bL  i_T^*(\tau_{<0}E) \to \bL  i_T^* E \to \bL  i_T^*(H^0E)$ we get
$$H^{i-1}\big(\bL  i_T^* (H^0E)\big)\ \  \to\ \  H^i\big(\bL
i_T^*(\tau_{<0}E)\big)\ \  \to\ \  H^i(\bL  i_T^* 
E) \ \ \to\ \  
H^i\big(\bL  i_T^*(H^0E)\big)\ \ \to \cdots$$ 
Note that, by right exactness of $\bL  i_T^*$ and left exactness of
$\bL  i_T^*[-1]$, we have that 
\begin{itemize}
\item $H^i\big(\bL  i_T^* (H^0E)\big) = 0$ for $i\neq 0,-1$,
\item $H^i \big(\bL  i_T^*(\tau_{<0}E)\big) =0 $ for $i>-1$.
\end{itemize}
also we have
\begin{itemize}
\item $H^i(\bL  i_T^* E)=0$ for $i\neq 0$ 
\end{itemize}
by assumption.

The long exact sequence implies that $H^i \big(\bL  i_T^*(\tau_{<0}E)\big) =0 $
for 
all $i$, so $\bL  i_T^*(\tau_{<0}E) =0$, and $E = H^0(E)$. \qed

\subsection{Invariance of the heart revisited}
We have proved earlier, that {\em assuming boundedness of $\cC$ with respect to
the standard $t$-structure, } the $t$-structure over a quasi-projective
variety is independent of the projective completion. We can now
remove the assumption:

\begin{proposition}\label{Prop:open-independent}
 Let $U\subset S_1$ be an open subset of a smooth
projective variety. Let $E \subset D(X\times U)$ and assume $E \in
\cC_{U\subset S_1}$. Then for every other open embedding $U\subset S_2$ in
a smooth projective variety, we have $E \in
\cC_{U\subset S_2}$.
\end{proposition}

{\bf Proof.} 
We apply induction on the dimension of $U$, the case of dimension 0
(or 1) 
being trivial.

{\sc Step 1: reductions.} 
By the sheaf property it suffices to prove the result for
a neighborhood of every point $s\in U$. 

Also, by Hironaka's theorem
(or any of its variants) we may choose a common resolution of
singularities of $S_1, S_2$ which is an isomorphism in a neighborhood
of any given $s\in U$, so it suffices to consider the case where
either $S_1 \das S_2$ or its inverse is a morphism; in either case
there are birational morphisms $S_i \to S'$ which are isomorphisms on
the same neighborhood. 

Finally we take a general divisor $D'$
on $S'$ through $s$ such that its pull-backs to both $S_1$ and $S_2$
are nonsingular. We denote its pullback to $U$ by $D$.  

{\sc Step 2: separating out the torsion.} 
Consider the maximal subobject $F \subset E$ with support in $X\times D$
 - see statement (2) in Proposition
\ref{Prop:max-torsion}. We have an exact sequence 
$$ 0 \to F \to E \to \bar E \to 0$$
  of objects in $\cC_{U\subset S_1}$. It suffices to show that $F$ and
$\bar E$ are in  $\cC_{U\subset S_2}$, since the latter category is closed
under extensions, being the heart of a $t$-structure.

{\sc Step 3: the non-torsion case.}  Let $D \subset D_1; \ D \subset
D_2$ be the
completions inside 
$S_1$ and $S_2$. Since $\bar E$ has no torsion in $D$, we have that 
$\bL i_D^* \bar E \in \cC_{D\subset D_1}$. The induction hypothesis
implies that $\bL i_D^* \bar E \in \cC_{D\subset D_2}$ as well. By the
open heart property (Proposition \ref{Prop:open-heart}) we have 
$\bar E_U' \in \cC_{U'\subset S_2}$ for some neighborhood $U'$ of
$D$, which is in particular a neighborhood of $s$ as required.

{\sc Step 4:  the torsion case.} Denote the defining section of $D$ by
$f \in H^0(U, \cO_U(D))$. Suppose $F$  is annihilated by $f^k$. We
have an exact sequence 
$$ 0 \to f\cdot F(-D) \to F \to F_0 \to 0$$
  of objects in $\cC_{U\subset S_1}$, where both $f\cdot F(-D)$ is
annihilated by $f^{k-1}$ and $F_0$ by $f$. It suffices to prove
the result for $f\cdot F(-D)$ and  $F_0$, and by induction on $k$ it
suffices to consider $k=1$.  So $F = i_{D\,*} F'$ where $F'\in
\cC_{D\subset D_1} = \cC_{D\subset D_2}$, so $F \in \cC_{U\subset
S_2}$ as needed.
\qed

\subsection{The generic flatness problem}

The following problem is fundamental:

\begin{problem} \label{Prob:generic-flatness}  Let $E \in \cC_S$. Is
there an open dense set 
  $U\subset S$ such that $E_U$ is $t$-flat over $U$?
\end{problem}

At this point we have only partial results.

\begin{lemma} If $S$ is a curve, and $E \in \cC_S$,  there is an
open dense set 
  $U\subset S$ such that $E_U$ is $t$-flat over $U$.
\end{lemma}

{\bf Proof.} Let $F\subset E$ be the maximal $S$-torsion subobject, and let $T$
be its support. Set  $U = S\setmin T$. Then $E_U$ is 
$S$-torsion-free and therefore $t$-flat.\qed

\begin{proposition}
 Let $E \in \cC_S$.  There is a
 dense subset  
  $Z\subset S$ such that $\bL i_s^* E\in \cC$ for all $s\in Z$.
\end{proposition}

{\bf Proof.} We apply induction on the dimension. First, we replace
$S$ by an open set where $E$ is $S$-torsion-free. Shrinking
$S$ further we may assume $S$ admits a smooth map $S \to
\PP^1$. The result holds for the fibers by induction, and therefore it
holds for $S$. \qed
 
According to \cite{Artin-Zhang}, the generic flatness problem has an
affirmative answer if the category $\cC_S$ can be defined and proved
noetherian for certain noetherian dedekind domains which are not
necessarily of finite type. At this point we know this to be true in
very few examples.

\section{Valuative criteria for semistable objects in
$\cP(1)$.}\label{Sec:valuative-P} 

\subsection{One parameter families for objects in $\cP(1)$}
Consider a numerical locally finite  stability condition  $(Z,\cP)$
with $\cP(\phi)
\subset D(X)$ and $Z:N(X) 
\to \CC$. Assume that its heart $\cC = \cP((0,1])$ is noetherian. We denote by
$(D(X)^{\le 0},D(X)^{\ge 0})$ the corresponding $t$-structure, where
$D(X)^{\le 0} = \cP((0,\infty))$ and $D(X)^{> 0} =
\cP((-\infty,0])$.

The main result of this section is the following:

\begin{theorem}\label{Th:valuative-P1} Let $S$ be a smooth curve and $U\subset
S$  a dense open subset.
 \begin{enumerate}
\item  Every family $F_U$ of objects in  $\cP(1)$ over $U$
 extends to a 
family $F$ of  objects in $\cP(1)$ over $S$.
\item
Let $F_1$ and $F_2$ be families of objects in $\cP(1)$ over $S$, and let
$\phi_U: (F_1)_U \to ( F_2)_U$ be an isomorphism. Then $\bL i^* F_1$ and
$\bL i^* F_2$ are S-equivalent.
\item Every family $F_U$ of objects in $\cP(1)$
over $U$ extends, after a suitable finite surjective base change $g:S' \to S$, to
a  
family $F$ of objects in  $\cP(1)$ over $S'$ with
polystable fibers in $S'\setmin g^{-1} U$.
\end{enumerate} 
\end{theorem}

Part (1) of the theorem follows from Proposition
\ref{Prop:extension-across-divisor}, in conjunction with 
the following lemma. We note that this lemma is the one point where our 
arguments restrict to  
$\cP(1)$ and do not extend to an arbitrary $\cP(t)$. 

\begin{lemma}\label{Lem:family-in-p(1)} Let $E$ be a family of objects in $\cC$
with connected base $S$ 
such that for
  some $s_0\in S$ we have $$\bL i_{s_0}^*E \in \cP(1).$$ 

Then for all $s\in S$ we have $\bL i_s^*E \in \cP(1).$
\end{lemma}
{\bf Proof.} The objects $\bL i_s^*E \in \cC$ are numerically
equivalent to each 
other.
Therefore  $Z(\bL i_s^*E) = Z(\bL i_{s_0}^*E) \in \bR _{<0}$. 
It follows that $\bL i_s^*E \in \cP(1).$ \qed

\subsection{Elementary modifications}

We proceed towards uniqueness in a standard manner. We start with elementary
modifications. For simplicity we work over an affine curve $S$, with
a fixed closed point $s$ having defining equation  $\pi$. We write
$i:\{s\} \to S$ and  $$U
= S \smallsetminus \{s\}.$$

Let $F\in \cC_S$ be a family of objects in $\cP(1)$. Consider an
exact sequence   
$$0 \to E \to  \bL i^*F \to Q \to 0.$$ Since $\bL i^*F \in \cP(1)$ we have that
$E,Q \in \cP(1)$. Since $F$ is $t$-flat, we have a
surjection $F \to i_*\bL i^*F$. Since $i_*$ is $t$-exact we have a
surjection $F \to i_* Q \to 0$. Let $G$ be its kernel.

\begin{definition}  We call $G$ the {\em elementary modification of $F$ at
$Q$.} 
\end{definition} 

\begin{lemma}
The object $G$ is a family of objects in $\cP(1)$. We have an exact sequence
$$0 \to Q \to \bL i^*G \to 
E 
\to 0,$$ in particular $\bL i^*G$ is S-equivalent to $\bL i^*F$. The
elementary modification of $G$ at $E$ is isomorphic to $F$, and 
the composition $F \to G \to F$ is the multiplication  $F \stackrel{\cdot \pi}{\to} F$ by the element $\pi$.
\end{lemma}

{\bf Proof.} Since $G$ is a subobject of $F$, the kernel of $\pi^n:G
\to G$ is a subobject of the kernel of  $\pi^n:F \to F$, which is 0. Thus $G$
is 
$t$-flat. By  Lemma \ref{Lem:family-in-p(1)} it is a family of objects in
$\cP(1)$.  

Since $i_*Q$ is annihilated by $\pi$, we have that $\pi: F \to F$ factors
through the inclusion $G \hookrightarrow F$. The fundamental
isomorphism $$\big(F/\pi F\big) \Big{/} \big(G/\pi F\big)\ \  \simeq\ \  F/G$$
 shows that $G/\pi F\simeq i_*E$,
and exchanging the roles of $F$ and $G$ we find that the kernel of $\bL i^*G
\to E$ is $Q$. \qed

\begin{lemma}
Let $\phi:F_A \to F_\Omega$ be a morphism of families of  objects of $\cP(1)$
giving an 
isomorphism 
$(F_A)_U \to ( F_\Omega)_U.$ Then there is a sequence of elementary
modifications $$F_A = F_0 \subset F_1 \subset \cdots \subset F_n = F_\Omega$$
whose composition is $\phi$. In particular $\bL i^*F_A$ and $\bL i^*F_\Omega$
are S-equivalent.
\end{lemma}

{\bf Proof.} The object $Q = F_\Omega/F_A$ has support at $\pi=0$, therefore
it is annihilated by $\pi^n$ for some positive $n$. Define $$K_i = \Ker
\, (\pi^{n-i}: Q \to Q),\ \  L_i= Q/K_i, $$ and let $F_i = Ker (F_\Omega
\to L_i)$. Then $F_{i+1} / F_{i}\simeq K_{i+1} /K_i$ is annihilated by
$\pi$. Thus $F_i$ is an elementary modification of $F_{i+1}$. The claim
follows from the previous lemma. \qed  

Using Lemma \ref{Lem:extend-morphism} we obtain:
\begin{lemma}
Let $F_A$ and $F_B$ be two families of objects in $\cP(1)$, and let $\phi:
F_{A|U} \to F_{B|U}$ be an isomorphism. Then $\bL i^*F_A$ and $\bL i^*F_B$ are
S-equivalent. 
\end{lemma}

\subsection{Polystable replacement}
For a polystable replacement we need a finite base-change. For this we
have
\begin{lemma}\label{Lem:flat-pullback}
Let $g: S'\to S$ be a finite flat morphism of nonsingular projective
varieties. Then 
$g_*$ and $ g^*$ are $t$-exact.
\end{lemma}

{\bf Proof.} Let $L$ be ample on $S$ and $L' = g^* L$. Denote by $p': X \times
S' 
\to X$ the projection. Consider an
object $E \in \cC_S$, and let $E' = g^* E$.  To show that $E' \in
\cC_{S'}$ we need to show that $\bR  p'_* (E' \otimes {L'}^n) \in \cC$ for large
$n$, which, by the projection formula  is the same as showing that $\bR  p_* (E
\otimes  g_*\cO_{S'}
\otimes L^n) \in \cC$, which is equivalent to $E \otimes g_*\cO_{S'}
\in \cC$, which follows from Proposition 
\ref{tensprop}, since $ g_*\cO_{S'}$ is a locally free sheaf.

Now let $E' \in \cC_{S'}$ and let us show that $ g_* E' \in \cC_S.$
We have by definition $\bR  p'_*E' \otimes {L'}^n \in \cC$ for large $n$,
so by the projection formula $\bR  p_*( g_* E'\otimes L^n) \in \cC$ for
large $n$, which means $ g_* E' \in \cC_S.$ \qed

\begin{lemma} 
Let $F$ be a family of objects in $P(1)$, with a sequence $$0 \to E \to
\bL  i^*F 
\to Q \to 0.$$ Consider its pullback via $g:S' \to S$ given by $\varpi^2 =
\pi$, and the incusion $i':\{s'\}\hookrightarrow S'$ of the point
$\varpi=0$. Let  
$H$ be the elementary modification of $ g^*F$ at $Q$. Then the exact
sequence  
$$0 \to Q \to \bL  {i'}^*H \to E \to 0$$ splits.
\end{lemma}

{\bf Proof.} Pulling back the exact sequence $$ 0 \to p^* Q \ \
\mathop\to\limits^\pi \ \  p^* Q
\to i_* Q\to 0 $$ and using Lemma \ref{Lem:flat-pullback} we get an exact
sequence  $$ 0 \to g^*p^* Q \ \ 
\mathop\to\limits^{\varpi^2} \ \  g^*p^* Q
\to g^*i_* Q\to 0. $$ 
We thus have an exact sequence 
$$0 \to i'_*Q \mathop\to\limits^\varpi g^* i_*Q \to i'_*Q \to 0, $$ and $H$ is the
pullback of $i'_*Q$ on the left along $g^*F \to  g^* i_*Q$. Thus we
have surjections $H \to 
i'_*Q$ and $H \to i'_*E$. 

The subobject $\varpi F \subset H$  surjects to $i'_*Q$
and maps to 0 in   $i'_*E$. 

Now consider $g^*G$, where $G$ is the elementary modification of $F$ at $Q$. We
have that $g^*G$ is canonically the elementary modification of $H$ at $i'_*Q$,
in particular it maps to 0 in  $i'_*Q$. At the same time it surjects to $g^*
i_* E$, in particular to  $i'_*E$. 

Altogether we have that $H$ surjects to   $i'_*Q\oplus i'_*E$, so ${i'}^*H \to
Q\oplus E$ is an isomorphism.

\qed

This completes the proof of Theorem \ref{Th:valuative-P1}.

\section{The noetherian property of discrete stability conditions}
\setcounter{subsubsection}{0}
In order for our results to be of some content, we need examples of
stability conditions with noetherian heart.

\begin{proposition}\label{Prop:discrete-is-noetherian}
Let $(Z,\cP)$ be a stability condition such that both $$Image (Z:K_0(X)
\to \CC)\subset \CC$$ and $$Image (Im Z:K_0(X)
\to \RR )\subset \RR $$ are discrete subgroups. Then the heart
$\cP((0,1])$ of the stability condition is noetherian.
\end{proposition}

\begin{corollary}
Let $T$ be a triangulated category of finite type with numerical rank
$2$. Then every numerical, locally finite stability condition on $T$
such that $\cP(0)\neq 0$ is noetherian. 

In particular, this holds in the following cases:
\begin{enumerate} 
\item $X$ is a smooth curve, and $T=D(X)$, and
\item $X$ a smooth projective variety, $\cC = Coh^{d,d-2}(X)$ is the
category of coherent sheaves supported in dimension $d$ modulo the Serre
subcategory of sheaves supported in dimension $d-2$, and $T =
D^b(\cC)$.
\end{enumerate}
\end{corollary}

{\bf Proof of corollary.} Let $(Z, \cP)$ be a numerical, locally
finite stability condition on $T$. If the image of $Z$ lies in a real
line, then local finiteness implies that $\cP(1) = \cP((0,1])$ is
noetherian. Otherwise the image of $Z$ is discrete, and since it
contains a  subgroup $\ZZ\cdot Z(\cP(0))$ of rank 1, the imaginary part
is discrete as well.   \qed

{\bf Proof of proposition.} We consider an increasing sequence 
$$ F_1 \subset F_2 \subset \cdots \subset E$$ of nonzero objects in the heart
$\cC=\cP((0,1])$.

\begin{lemma} Suppose $E \in \cP(1)$. Then $F_i\in \cP(1)$ and the
sequence stabilizes.
\end{lemma}

{\bf Proof.} We  have $$Im(Z(F_i)) = - Im Z(E/F_i)$$ and both are
nonnegative, therefore $Im(Z(F_i)) =0 $ and thus $F_i \in \cP(1)$. 
Now $Z(F_i)$ and $Z(F_{i+1}/F_i)$ are nonpositive, therefore $Z(F_i)$
is nonincreasing. For the same reason $Z(F_i) > Z(E)$, and by discreteness
$Z(F_i)$ stabilizes. But at that point $Z(F_{i+1}/F_i) = 0$ and thus
$F_{i+1}/F_i = 0$, so the sequence stabilizes.
\qed

\begin{lemma}
Let $E_0\subset E$ with $E_0\in \cP(1)$ be the (possibly zero)
Harder--Narasimhan subobject of phase 1, and let $\phi_{\max}(E/E_0) =
a <1$ be the phase of the next Harder--Narasimhan piece. Let $T_i 
\subset F_i$ be the maximal subobject with Harder--Narasimhan
constituents of phases $>a$. 

Then $T_i\subset T_{i+1} \subset E_0$ form an increasing sequence and
this sequence stabilizes. 
\end{lemma} 

{\bf Proof.} The morphism $T_i \to E/E_0$ vanishes since
$\phi_{\min}(T_i) > \phi_{\max}(E/E_0) = a$. Therefore $T_i
\hookrightarrow E$ factors through $E_0$. Similarly $T_i
\hookrightarrow F_{i+1}$ factors through $T_{i+1}$. This sequence
stabilizes by the previous lemma.
\qed

{\bf Proof of proposition.} By the previous lemmas we may assume that
$E \not\in \cP(1)$ and $T_i=T$ are constant. Consider the increasing
sequence of subobjects 
$F_i/T \subset E/T$. By the discreteness assumption we have that the
nondecreasing sequence  $Im\
Z(F_i/T)$ stabilizes, at some $b>0$. So we may assume it is
constant. Then $Re\ 
Z(F_i/T)$ is nonincreasing, but $\phi_{\max}(F_i/T) \leq a$. Therefore
$Z(F_i/T)$ stabilizes, and as before this implies that $F_i/T$
stabilizes and thus $F_i$ does.
\qed

Recall that to every connected component $\Sigma$ in the space of numerical
stability conditions on $T$ Bridgeland associates a subspace
$V(\Sigma)\subset(\cN(T)\otimes\CC)^*$, such that the map sending a stability
to the corresponding central charge induces a local homeomorphism
$\Sigma\to V(\Sigma)$. It seems to be unknown whether the subspace $V(\Sigma)$
is always defined over $\QQ$. This is the case for all
the components explicitly described in \cite{Bridgeland} and \cite{Bridgeland2}
for the cases of curves and K3 surfaces.

\begin{corollary}\label{rat-cor} 
Let $\Sigma$ be a connected component in the space of numerical
stability conditions for which $V(\Sigma)$ is defined over $\QQ$.
Then the set of stability conditions satisfying the assumptions of Proposition
\ref{Prop:discrete-is-noetherian} is dense in $\Sigma$.
\end{corollary}

{\bf Proof.} Indeed, if $V(\Sigma)$ is defined over $\QQ$ then the set
of homomorphisms $Z:\cN(T)\to\CC$ with the image contained in $\QQ+i\QQ$ is
dense in $V(\Sigma)$.
\qed 

\section{Questions} \label{Sec:questions}

\subsection{Towards a moduli space of semistable objects}
In order to construct proper moduli spaces of S-equivalence classes of
objects in $\cP(1)$ with fixed numerical class, two problems remain:
\begin{enumerate}
\item The generic flatness problem (Problem
  \ref{Prob:generic-flatness}).  
\item Boundedness: fixing a  class $c\in H^*(X,\QQ)$ there
  should be a scheme of finite type $S$ and a family $E \in 
  D(X\times S)$ including all objects in $\cP(1)$ having Chern character $c$. 
\end{enumerate}
Generic flatness would imply openness of versality using the open heart
property (Proposition \ref{Prop:open-heart}). If this holds, then
by the work of Lieblich \cite{Lieblich} (generalizing the work
of Inaba \cite{Inaba}) these objects form an Artin algebraic stack.
The moduli space of S-equivalence classes of
objects in $\cP(1)$ can be pieced together from this stack -
this will be pursued elsewhere.

 Boundedness is necessary for proving
properness, and we have not begun pursuing the various approaches to
address this. It is expected to be a difficult problem. It was, however, noted by Tom Bridgeland that if the class of objects of bounded mass and amplitude subject to one stability condition is a bounded family, then the same is true for all stability conditions in the same connected component of the space of stability conditions. It is expected that this kind of boundedness should be easier to show near the "large volume limit".

Individual moduli spaces may be constructed directly without relying on such results. It has been claimed (without proof, as far as we know) that Bridgeland's moduli space of perverse point sheaves is an instance of this phenomenon. We have been informed that D. Arcara and A. Bertram are in the process of working out a new non-trivial example on a K3 surface which is related by a Mukai flop to a well known moduli space of curvilinear torsion sheaves.

\subsection{Alternative approaches}
Our construction of the sheaf of
$t$-structures is rather roundabout, and given the rather simple
characterization of the resulting $t$-structure one may wonder if
there is no direct way to this construction. It is conceivable that a
further study of the quasi-coherent picture as in \cite{LST}, may lead to
a direct construction, which may be more general. 

Given an abelian category with a number of assumptions, Artin and
Zhang \cite{Artin-Zhang} study various aspects of moduli problems for
objects in this 
category, in particular $Quot$ functors. It is interesting to see to
what extent our method can be combined with theirs (or maybe even overridden
by their approach).  
\subsection{Removal of superfluous assumptions}
Various aspects of our discussion requires assumptions which might not
be necessary. 
\begin{enumerate}
\item It would be nice to extend the results to $\cP(t)$ for $t\neq 1$.
\item The characterization of the $t$-structure over a
quasi-projective variety should be general, and should require neither
boundedness
of $\cC$ with respect to the standard $t$-structure
nor resolution of singularities.
\item it would be interesting to see to what extent the noetherian
assumption may be weakened. After all, Bridgeland's stability
conditions are only assumed locally finite in general.

One needs to be a bit careful: some condition is necessary for any sort of extension result, as shown by the following example - suggested by both the referee and by Bridgeland, who attributed it to R. Thomas:

Consider $X=\PP^1, S= \bbA^1$. Let $(T,F)$ be the tilting pair on $\PP^1$ with $T=$ torsion supported at $0\in \PP^1$. The heart $\cC$ of the $t$-structure given by tilting this pair, with an appropriate grading, consists of complexes $\cF\in D(\PP^1)$ with $H^0(\cF)\in F, H^1(\cF)\in T$. This is clearly not noetherian: $\cO \to \cO(1)\to \cO(2)\to\cdots$ is an infinite sequence of $\cC$-epimorphisms, given by the defining section of the point $0\in \PP^1$.  

Take the structure sheaf $\cG^0 = \cO_{\Delta^0}$ of the diagonal $\Delta^0\subset X \times (S\smallsetminus \{0\})$. Its fibers are in $\cC$, but what could be an extension $\cG\in D(X\times S)$ with fibers in $\cC$? any extension is supported on the union of the diagonal with $X \times \{0\}$, it is easy to see that whatever has support on the generic point $\Spec K(X) \times \{0\}$ cannot give something in $\cC$, and one is left with an object supported on the diagonal, whose standard $H^0$ has torsion at 0, thus it is not in $\cC$. In other words - an extension in $\cC$ does not exist.

\item It would be interesting to see to what extent the construction of the $t$-structure can be extended to singular base schemes. Our construction breaks down exactly where we describe the essential image under push-forward of a closed subvariety, but we do not have an example where a truncation functor fails to exist.
\end{enumerate}  

Techniques such as in
\cite{Bondal-vandenBergh,Kashiwara,Yekutieli-Zhang} may prove useful
in pursuing these directions.

\end{document}